\documentclass[10pt]{article}
\usepackage[margin=2.5cm]{geometry}
\usepackage[utf8]{inputenc}
\usepackage{bussproofs}
\usepackage{calc}
\usepackage{amsmath}
\usepackage{amsthm}
\usepackage{amssymb}
\usepackage{stmaryrd}
\usepackage{setspace}
\usepackage{MnSymbol}
\usepackage{cmll}
\usepackage{longtable}
\theoremstyle{definition}
\newtheorem{thm}{Theorem}
\newtheorem{cnv}[thm]{Convention}
\newtheorem{defn}[thm]{Definition}
\newtheorem{crl}[thm]{Corollary}
\newtheorem{cnj}[thm]{Conjecture}
\newtheorem{prp}[thm]{Proposition}

\theoremstyle{remark}
\usepackage{setspace}
\title{Denotational semantics for languages of epistemic grounding\\based on Prawitz's theory of grounds}
\author{Antonio Piccolomini d'Aragona\\\begin{small}Aix Marseille Univ, CNRS, Centre Gilles Gaston Granger, Aix-en-Provence, France\end{small}\\\texttt{antonio.piccolomini-d-aragona@univ-amu.fr}}
\date{}

\begin{document}

\maketitle

\begin{abstract}
\noindent We outline a class of term-languages for epistemic grounding inspired by Prawitz's theory of grounds. We show how denotation functions can be defined over these languages, relating terms to proof-objects built up of constructive functions. We discuss certain properties that the languages may enjoy both individually (canonical closure and universal denotation) and with respect to their expansions (primitive/non-primitive and conservative/non-conservative expansions). Finally, we provide a ground-theoretic version of Prawitz's completeness conjecture, and adapt to our framework a refutation of this conjecture due to Piecha and Schroeder-Heister.
\end{abstract}

\begin{small}
\end{small}

\section{Introduction}

Prawitz's \emph{theory of grounds} (ToG) is a \emph{proof-theoretic semantics} centered on the notion of \emph{ground} [Prawitz 2009, 2012a, 2013, 2015, 2019 - the expression ‘‘proof-theoretic semantics" is due to Schroeder-Heister 1991; for a general introduction see Francez 2015, Schroeder-Heister 2018]. Its main task is that of explaining why and how valid inferences and proofs have the power to epistemically compel us to accept their conclusions, if we have accepted their premises.

ToG shares many aspects with Prawitz's previous semantics of valid arguments and proofs [Prawitz 1973, 1977], as well as with some of Prawitz's earliest papers on constructive semantics [Prawitz 1971a]. However, it also introduces a number of structural and formal novelties that, with respect to the task of accounting for epistemic compulsion, permit Prawitz to overcome some problems his earlier approaches suffered from. The latter are unsatisfactory principally because of the explanatory order of the intertwinement of the notions of valid inference and proof: valid inferences are defined in terms of preservation of validity throughout the global proof-structures where they occur. For example, in Prawitz's semantics of valid arguments, the validity of a conjunction-elimination

\begin{prooftree}
\AxiomC{$\alpha_1 \wedge \alpha_2$}
\UnaryInfC{$\alpha_i$}
\end{prooftree}
($i = 1, 2$) is explained by saying that, given a closed valid argument $\Delta$ for $\alpha_1 \wedge \alpha_2$, the argument

\begin{prooftree}
\AxiomC{$\Delta$}
\noLine
\UnaryInfC{$\alpha_1 \wedge \alpha_2$}
\UnaryInfC{$\alpha_i$}
\end{prooftree}
is valid. This is because, in this semantics, a closed valid argument reduces by definition to a closed argument for the same conclusion, whose last step is an introduction and whose immediate sub-arguments are valid. So we have

\begin{prooftree}
\AxiomC{$\Delta_1$}
\noLine
\UnaryInfC{$\alpha_1$}
\AxiomC{$\Delta_2$}
\noLine
\UnaryInfC{$\alpha_2$}
\BinaryInfC{$\alpha_1 \wedge \alpha_2$}
\UnaryInfC{$\alpha_i$}
\end{prooftree}
At this point, conjunction elimination can be ‘‘justified" by means of a rule for removing so-called maximal conjunctions, i.e. occurrences of the conjunction which are both conclusions of an introduction and (major) premises of an elimination:

\begin{prooftree}
\AxiomC{$\Delta_1$}
\noLine
\UnaryInfC{$\alpha_1$}
\AxiomC{$\Delta_2$}
\noLine
\UnaryInfC{$\alpha_2$}
\BinaryInfC{$\alpha_1 \wedge \alpha_2$}
\UnaryInfC{$\alpha_i$}
\AxiomC{reduces to}
\noLine
\UnaryInfC{}
\noLine
\UnaryInfC{}
\noLine
\UnaryInfC{}
\AxiomC{$\Delta_i$}
\noLine
\UnaryInfC{$\alpha_i$}
\noLine
\UnaryInfC{}
\noLine
\UnaryInfC{}
\noLine
\TrinaryInfC{}
\end{prooftree}
So, the validity of our conjunction elimination ultimately depends on the validity of $\Delta$. And clearly, this is problematic if we also want to explain the epistemic power of $\Delta$ in terms of the epistemic power of the valid inferences which $\Delta$ is built up of - including our conjunction elimination. In ToG, therefore, the explanatory order is reversed. Inferential validity is no longer defined globally but locally, through a non-inferential notion of ground. Grounds are meant to reify what one is in possession of when in a state of evidence for judgments or assertions. This also requires to conceive of inferences, not as simple transitions from certain premises to a certain conclusion, but more strongly as applications of constructive operations on grounds that one considers oneself to have for the premises, with the aim of obtaining a ground for the conclusion. An inference is valid if the operation on grounds it involves actually produces grounds for the conclusion from grounds for the premises. Thus, grounds will be built up of the operations applied in inferences, and they will thereby remind us of proof-objects of propositions in Martin-L\"{o}f's intuitionistic type theory [Martin-L\"{o}f 1984]. Going back to our example above, the validity of

\begin{prooftree}
\AxiomC{$\alpha_1 \wedge \alpha_2$}
\UnaryInfC{$\alpha_i$}
\end{prooftree}
is now explained by defining a function, say $f_{\wedge, i}$, such that, for every ground $g$ for asserting $\alpha_1 \wedge \alpha_2$, $f_{\wedge, i}(g)$ is a ground for asserting $\alpha_i$. This is done through an equation that renders $f_{\wedge, i}$ total and constructive, i.e. given that a ground for asserting $\alpha_1 \wedge \alpha_2$ involves by definition a pair whose $i$-th element is a ground $g_i$ for asserting $\alpha_i$,

\begin{center}
    $f_{\wedge, i}(g) = g_i$.
\end{center}
The basic idea is that the inference is compelling since a performance of it corresponds to an aware application of $f_{\wedge, i}$ to $g$; by applying $f_{\wedge, i}$ to $g$, we ‘‘see" that a ground for asserting $\alpha_i$ obtains. Proofs can be now defined as chains of valid inferences, and their epistemic power can be explained, in a non-circular way, in terms of the epistemic power of the valid inferences of which they are made up [further details about the reasons that led Prawitz from his previous approaches to ToG can be found d'Aragona 2019a, 2019b, Tranchini 2014a, Usberti 2015, 2019].\footnote{Two remarks. First, one may very well say that a ground is nothing but a proof, so our explanation would remain circular. This is a good objection, but the answer is that in ToG Prawitz endorses a proof-objects/proof-acts distinction. While grounds are proof-objects, proofs are proof-acts. What we have to explain is the epistemic compulsion exerted by acts, and this could not be done in Prawitz's previous framework \emph{exactly} because proofs and valid arguments were there understood \emph{both} as proof-objects \emph{and} as proof-acts [for further details, see d'Aragona 2021b]. Second, the idea that one ‘‘sees" that a ground for the conclusion of a valid inference is obtained by applying the relevant operation is just a \emph{desideratum}. It is not clear whether and how this \emph{desideratum} can be actually achieved. Of course, it all depends on how one defines operation on grounds, and in particular, on whether this definition makes it in some sense recognizable that an operation on grounds yields grounds for asserting a certain formula when applied to grounds for asserting given formulas. ToG raises a recognizability problem which also affects Prawitz's previous semantics in terms of proofs and valid arguments [see mainly Prawitz 2015] and which is sometimes referred to BHK-semantics too [see Díez 2000]. We discuss this point in Section 9.3 below.}
 
Prawitz [Prawitz 2015] also envisages the development of formal languages of grounding, whose terms, so to say, encode or describe deductive activity, and denote the grounds it produces. These languages, Prawitz says, are to be a sort of extended $\lambda$-languages, except that they must be understood as indefinitely open to the addition of new linguistic resources for new operations. In fact, because of G\"{o}del's incompleteness theorems no recursively axiomatisable language of grounding over first-order arithmetic will be able to account for all the possible first-order arithmetic grounds; an undecidable sentence will have a ground not expressible in this language, but in a suitable expansion of it.

ToG offers an interesting and original insight into deduction and its principles, merging important suggestions in the constructivist field such as BHK semantics [Troelstra \& Van Dalen 1988], the \emph{formulas-as-types conception} of the Curry-Howard isomorphism [Howard 1980], Dummett's investigations into theory of meaning [Dummett 1991, 1993a, 1993b], proof-theory originating with Gentzen [Gentzen 1934 - 1935] and Prawitz's own normalisation theory [Prawitz 2006]. In addition, it shows that grounding - a widespread topic in contemporary logic [see Poggiolesi 2016, 2020 for an overview] - can be addressed from an epistemic standpoint alien to most of the approaches. However, we remark that in Prawitz the word ‘‘grounding" has not the same meaning as in the metaphysical literature, for Prawitz does not focus on propositions or sentences grounding another proposition or sentence, but rather on \emph{how} such links obtain.

Despite these indications, however precious, the development of formal languages of grounding is only suggested, but not accomplished in Prawitz's ground-theoretic writings. Here, we work upon these suggestions. First, we define general notions of language of grounding and of expansion of a language of grounding. Second, we introduce denotation functions that associate terms of languages of grounding to grounds and operations on grounds. Finally, we investigate certain properties that denotation functions induce on languages and their expansions. Our languages contain terms only, but in the concluding remarks we suggest a possible enrichment of them through the addition of formulas, and a parallel enrichment of our semantics through the formulation of models based on truth- or proof-conditions. As a step preliminary to all this, however, we describe the ‘‘universe" of grounds and operations on grounds that the terms have to denote, as well as some standard concepts in Prawitz's semantics like that of atomic base. The proposed setup has many links with Prawitz's semantics of valid arguments, or with other approaches based on the Curry-Howard correspondence - e.g. Martin-L\"{o}f's type-theory [see for example Martin-L\"{o}f 1984]. But there will also be deep differences, due to will of complying with Prawitz's philosophical purposes in ToG.\footnote{Prawitz's proposal has been taken up also by [d'Aragona 2021a], where a concrete language of grounding is exemplified and used to indicate how a language of grounding is to be made in general. [d'Aragona 2021a] suggests that languages of grounding could be equipped with denotation functions associating their terms to grounds and operations on grounds, and that a hierarchy of these languages may be created by introducing an extension-relation among them and denotation functions defined over them. However, this strategy is only hinted at in [d'Aragona 2021a], whose real aim is rather that of introducing formal systems through which properties of terms of the underlying languages of grounding can be proved. For this reason, the languages of grounding introduced in [d'Aragona 2021a] do not contain terms only, but also formulas expressing such things as ‘‘the given term denotes a ground for the given formula" or ‘‘the given terms are equivalent". The hierarchy of languages is hence understood as a hierarchy of languages and systems over these languages.}

\section{General overview}

We start with an informal overview of some general concepts underlying our formal proposal. 

\subsection{Forms of judgments and assertions}

As said, Prawitz's grounds reify evidence for judgments and assertions. Thus, we must first of all say what judgments and assertions are. Our line of thought is essentially inspired by Prawitz's own standpoint [see mainly Prawitz 2015].

A judgment can be understood as a mental act where one claims a proposition to be true. And just like sentences can be looked at as the linguistic counterpart of propositions, assertions can be taken to be the linguistic counterpart of judgments. An assertion is hence a linguistic act where one claims a sentence to be true. This is not meant as implying that judgments or assertions are of the form ‘‘$\dots$ is true", where the dots are replaced by (a name of) a proposition or sentence. The claim of truth is implicit; judgments or assertions just employ propositions or sentences in a certain mood, or with a certain force.

Philosophically, the notions of judgment and assertion can of course be said to be very different. Nonetheless, Prawitz declares that it is ground-theoretically irrelevant whether an inference takes place at a mental or else at a linguistic level [see for example Prawitz 2015]. Here, we shall endorse this line of thought, and consider judgments and assertions as interchangeable. In particular, we shall limit ourselves to speaking of assertions.

Assertions of sentences are called by Prawitz \emph{categorical}. We indicate them through the fregean (meta-linguistic) turnstile, i.e. $\vdash \alpha$, for $\alpha$ closed. Prawitz also considers other kinds of assertions: \emph{general} assertions, i.e. assertions involving open formulas, noted
\begin{center}$\vdash \alpha(x_1, ..., x_n)$;
\end{center}
\emph{hypothetical} assertions, i.e. assertions involving a sentence depending on other sentences, noted
\begin{center}$\alpha_1, ..., \alpha_n \vdash \beta$;
\end{center}
\emph{general-hypothetical} assertions, i.e. assertions involving a possibly open formula that depends on other possibly open formulas, noted
\begin{center}$\alpha_1(\underline{x}_1), ..., \alpha_n(\underline{x}_n) \vdash \beta(\underline{x})$,
\end{center}
with $\underline{x}_i$ and $\underline{x}$ sequences of the variables occurring free in $\alpha_i$ and $\beta$ respectively ($i \leq n$), where at least one of these sequences is not empty.

We call of \emph{first-level} the assertions illustrated so far. Assertions of \emph{second-level} can again be \emph{hypothetical} or \emph{general-hypothetical}, and involve a possibly open formula, a hypothetical or a general-hypothetical assertion of first-level depending on at least one hypothetical or general-hypothetical assertion of first-level - and maybe other possibly open formulas. The notation is in this case
\begin{center}
$\tau_1, ..., \tau_n \vdash \tau_{n + 1}$
\end{center}
where $\tau_i$ and $\tau_{n + 1}$ are possibly open formulas, or general or general-hypothetical assertions of first-level, and where we have the restrictions: there is at least one $\tau_i$ hypothetical or general-hypothetical of first-level; if $\tau_{n + 1}$ is hypothetical or general-hypothetical of first-level, what is on the left of its turnstile is contained in or equal to the union set of what is on the left of the turnstile of the hypothetical or general-hypothetical $\tau_i$-s of first-level ($i \leq n$).

\subsection{Bases and (operational) types}

The objects of ToG, Prawitz says, are \emph{typed}, i.e. grounds and operations on grounds have as \emph{type} the sentences involved in the assertions they either are a ground for, or the grounds for which they take and produce as values [see mainly Prawitz 2015]. Since we aim at developing formal languages of grounding whose terms denote such objects, these terms will have to be typed as well. To begin with, therefore, we need a background language from which types are drawn. Following Prawitz, we take this background language to be a first-order logical language, so that we leave aside the question about whether our proposal - and Prawitz's ToG - can be extended also to higher-order logics.

With this language established, we have to know what its formulas mean. In fact, since we are speaking of evidence, and since the latter must in the end depend on what the sentences involved in assertions mean, grounds and operations on grounds will be such \emph{relatively} to the meaning of their type - and they may be qualified as logical if they are so \emph{independently} of the meaning of the non-logical signs. While the meaning of the logical constants is fixed in terms of appropriate semantic clauses, the meaning of the non-logical signs will be fixed in terms of atomic systems that show how these signs behave in deduction. Reference to atomic bases is customary in Prawitz-style approaches, being for example extensively used also for valid arguments. However, the idea that atomic systems determine meaning is anything but straightforward. Complex notions are at play here, like that of a rule ‘‘concerning" a symbol, and more rules ‘‘concerning" one or more symbols according to some order. Unfortunately, a detailed treatment of these difficult issues would lead us too far; apart from some examples, we will therefore leave it aside here, and take the notion of rules ‘‘concerning" symbols as somewhat primitive [the interested reader can refer to Cozzo 1994, where a theory of meaning is proposed, meant to answer these questions].

The idea that grounds are typed on formulas of a first-order language clearly implies the full adoption of the \emph{formulas-as-types} conception of the Curry-Howard isomorphism [Howard 1980] - something that already Prawitz does for his own ToG [see mainly Prawitz 2015]. But since we want to speak also of operations from grounds to grounds, we have to speak more generally of operational types, i.e. types with domain, co-domain, and possibly free variables. As a matter of fact, this amounts instead to requiring a kind of dependent type theory in the style of Martin-Löf [see for example Martin-Löf 1984].

\subsection{Summary scheme}

The first step of our proposal will be that of defining a domain of grounds and operations on grounds. These objects are classified according to the ‘‘complexity" of their type. The basic notion is that of ground for a categorical assertion.

Categorical assertions require \emph{conclusive} evidence. In Frege's terminology, evidence must in this case be reified by a \emph{saturated} entity, that we shall call a (\emph{closed}) \emph{ground}. So, a ground for an assertion $\vdash \alpha$, for $\alpha$ closed with main logical sign $\star$, will have type $\alpha$, and will be fixed through a primitive, meaning-constitutive operation on grounds $\star I$, via clauses that run by induction on the complexity of $\alpha$.

We can then move on to (\emph{open}) \emph{grounds} for general, hypothetical and general-hypothetical assertions of first-level. The general idea here is that of an \emph{unsaturated} entity, i.e. of an effective method that applies, first, to individuals corresponding to the individual variables occurring free in the assertion (if any), and then, to grounds for the categorical assertions obtained from the assertions on the left of the turnstile (if any) by replacing the free individual variables with names of the individuals. The method produces a ground for the categorical assertion obtained from the assertion on the the right of the turnstile by replacing free individual variables with names of the individuals. We shall speak in this case of an operation on grounds of \emph{first-level}, whose operational type has the assertions on the left of the turnstile as domain, and the assertion on the right of the turnstile as co-domain.

Finally, we can deal with operations on grounds of \emph{second-level}, also called (\emph{open}) \emph{grounds} for the corresponding hypothetical or general-hypothetical assertions. These are again \emph{unsaturated} entities, except that the methods required in this case take as values, not only grounds for categorical assertions, but at least one operation on grounds of first-level - possibly plus individuals on a reference domain. And they generate, either grounds for categorical assertions, or an operation on grounds of first-level.

The second step of our proposal is the introduction of a hierarchic class of formal languages of grounding. As already said, we cannot limit ourselves to ‘‘recursive" languages of grounding, i.e. languages for denoting recursive sets of grounds and operations on grounds. If this were the case, because of G\"{o}del's incompleteness we would be unable to account for all the grounds on atomic bases as rich as a base for first-order arithmetic. On the other hand, given a ‘‘recursive" language of grounding on an atomic base for, say, first-order arithmetic, any ground or operation on grounds not captured through this language can be described by suitably expanding the language - for example, by adding to the base some ground-theoretic reflection principle.\footnote{Observe that, strictly speaking, we have so far referred to what in Section 8 we will call \emph{non-conservative} expansions, namely, expansions of languages of grounding obtained by extending the atomic base. We have argued that G\"{o}del's incompleteness results imply that we \emph{need} expansions of this kind. However, there are reasons for maintaining that we also need what, again in Section 8, we will call \emph{conservative} expansions, namely, expansions obtained by adding symbols for new operations on grounds \emph{on the same base}. In fact, a language of grounding may not be ‘‘(strongly) complete" with respect to the grounds and operations on grounds on its atomic base. In other words, we have no guarantee that any ground or operation on grounds on this base can be described as a combination of the grounds and operations on grounds denoted by the terms of our language nor, more weakly, that for any (operational) type inhabited in the base, there is a ground or operation on grounds denoted by some term of our language that inhabits that (operational) type.} The expansion relation will not in general hold between any two languages, even if the smaller language is defined over a sub-base of the atomic base of the bigger language - e.g., typically, if in the base there are no grounds for $\vdash \alpha$, there is an operation on grounds on this base from grounds for $\vdash \alpha$ to grounds for $\vdash \beta$ for any $\beta$, but this will cease to hold as soon as we extend the base so as to have grounds for $\vdash \alpha$. Alternatively, one may define languages of grounding which be ‘‘maximal" on an atomic base, namely, (possibly non-recursive) languages whose terms denote all the grounds and operations on grounds on the base. An approach of this kind is put forward in [Prawitz 2014]. Here we prefer to stick to the hierarchic standpoint. It will permit us to deal with languages which, as we shall see below, can be understood as recursive formal systems where non-introductory inferences are equipped with justification procedures. The properties of our languages can be therefore seen as the properties of ground-theoretically sound systems ordered according to, and progressively improved via some underlying semantic principles.

Over our languages, we define a denotation function in such a way that the typed terms of the language be ‘‘names" of the ‘‘inhabitants" of the domain of grounds or operations on grounds. The pattern is as follows. Closed terms of closed type $\alpha$ denote grounds for categorical $\vdash \alpha$; open terms of open type $\alpha(x_1, ..., x_n)$ with $n \neq 0$ individual variables occurring free denote grounds for general $\vdash \alpha(x_1, ..., x_n)$; open terms with $n \neq 0$ typed-variables $\xi^{\alpha_1}, ..., \xi^{\alpha_n}$ occurring free, for closed $\alpha_i$ ($i \leq n$) and of closed type $\beta$, denote grounds for hypothetical $\alpha_1, ..., \alpha_n \vdash \beta$; open terms with $n \neq 0$ individual variables occurring free, and $m \neq 0$ typed-variables $\xi^{\alpha_1(\underline{x}_1)}, ..., \xi^{\alpha_m(\underline{x}_m)}$ occurring free, for closed or open $\alpha_i(\underline{x}_i)$ ($i \leq n$) and of closed or open type $\beta(\underline{x})$, denote grounds for general-hypothetical $\alpha_1(\underline{x}_1), ..., \alpha_m(\underline{x}_m) \vdash \beta(\underline{x})$. As for operations on grounds of second-level, they will be denoted by no term - this is no limitation, for one can attribute a denotation also to the operational symbols that terms are built up of [a similar approach is in d'Aragona 2018].

\subsection{Relations with other approaches}

The approach outlined in the previous section requires a simultaneously recursive definition of the notions of ground and operation on grounds. This occurs already in Prawitz's semantics of valid arguments. The main idea is that the open cases be explained through their closed instances - what Schroeder-Heister [Schroeder-Heister 2008, 2012] and Do\v{s}en [Do\v{s}en 2015] have called the \emph{dogma of the primacy of the categorical over the hypothetical}.

The (operations on) grounds and the terms that we deal with below can be seen as ‘‘functional" translations of the various kinds of Prawitz's valid arguments: closed grounds correspond to closed valid arguments in full-evaluated form; closed terms correspond to closed valid arguments; open terms, and the operations they denote, correspond to open valid arguments. More generally, operations on grounds can be considered as a ‘‘functional" version of valid generalized inference rules up to order $2$ as defined by Schroeder-Heister [Schroeder-Heister 1984a, 1984b]. It is hence not surprising that the type structure, the objects and the terms of our framework will look like the type structure, the objects and the terms of other approaches based on the Curry-Howard correspondence. In particular, the type structure comes very close to that of Martin-L\"{o}f's dependent type theory while - as seen from the meaning explanation of the previous section - objects and terms are nothing but proof-objects for propositions in intuitionistic type theory [Martin-L\"{o}f 1984].

Such similarities are acknowledged, and even required by Prawitz himself [see mainly Prawitz 2015]. In fact, they are not limited to the formal level of types, objects and terms, but also concern more philosophical aspects. In ToG, for example, Prawitz introduces a distinction between proof-objects (grounds and operations on grounds) and proof-acts (proofs, i.e. chains of valid inferences, i.e. chains of applications of operations on grounds). As known, this distinction plays a pivotal role in Martin-L\"{o}f's intuitionistic type theory, where it is understood more in particular as a distinction between proofs of propositions and demonstrations of judgments [see for example Martin-L\"{o}f 1984; for a thorough and seminal analysis of the proof-objects/proof-acts distinction, see Sundholm 1998]. In spite of this and of other similarities, though, Prawitz's and Martin-L\"{o}f's views still differ in many fundamental respects. At variance with Martin-L\"{o}f (and Sundholm), Prawitz does not consider proof-objects as devoid of any epistemic content [see Prawitz 2012b]. For him, a ground is \emph{not only} a truth-maker of $\alpha$, but also the denotation of a term that \emph{encodes} a proof of $\alpha$, and hence of something that justifies the assertion of $\alpha$. In turn, the latter is a mere speech-act, not an expression of the form ‘‘$\alpha$ is true", which Prawitz would consider as belonging to a meta-level [see Prawitz 2015].

The main task of ToG is that of defining inferential validity, following the idea that an inference is valid when an operation exists from grounds for the premises into grounds for the conclusion. Thus, even if ToG \emph{could} be developed as a theory where one speaks of how given terms denote grounds, it should be primarily understood as a theory about \emph{what} grounds and operations on grounds are. And since we have knowledge of grounds and of operations on grounds only through adequate linguistic descriptions of them, ToG should also be primarily understood as a theory \emph{of} the terms that denote grounds, and so of the languages where these terms can be built.

It is for this reason that, in spite of the strong similarities between the type structure and the terms that we deal with below, on the one hand, and the type structure and the proofs of propositions in Martin-L\"{o}f's type theory on the other, the line of thought we pursue is quite different in spirit from that of a Martin-L\"{o}f-based approach. We do not develop theories for proving that terms denote grounds, but a framework for the languages over which such theories can be built. We do not have rules that produce types, but assume types as given. We do not take the meaning of the types as ‘‘embedded" into the theory, but set it at a meta-level, via objects called grounds as well as through denotation functions that associate terms to these objects. In a sense, our proposal comes rather closer to a theory of \emph{justified deductive systems} as defined in [Prawitz 1973] - as we shall see, a language of grounding over which a denotation function is defined can be understood as a Curry-Howard translation of a recursive system with rules equipped of appropriate justification procedures.

\section{Preliminary definitions: bases and types}

In this section we introduce the notions of background language - i.e. a first-order logical language - and of atomic base. We remark that we shall limit ourselves to first-order languages - as Prawitz himself usually does - thus leaving aside the question whether our proposal can be extended to higher-order languages. Clearly, this is not to deny that such an extension would be of interest.

\subsection{First-order logical languages}

\begin{defn}
A \emph{first-order language} $L$ is defined as usual by:

\begin{itemize}
    \item an alphabet $\texttt{Al}_L$ with denumerably many individual variables, a (possibly empty) set of individual constants, a set of relational constants, an atomic logical constant $\bot$ for the absurd, logical constants $\wedge$, $\vee$, $\rightarrow$, $\forall$, $\exists$, and auxiliary symbols;
    \item a set $\texttt{TERM}_L$ of terms specified by induction;
    \item a set $\texttt{FORM}_L$ of formulas specified by induction, with a subset $\texttt{ATOM}_L$ of atomic formulas specified in a standard way.
\end{itemize}
The context will clarify whether meta-variables - possibly with indexes - indicate individual variables, individual constants, arbitrary terms, formulas and sets of formulas. Negation is defined as

\begin{center}
$\neg \alpha \stackrel{def}{=} \alpha \rightarrow \bot$
\end{center}
Over $L$, we assume as defined in a standard way:

\begin{itemize}
    \item free variables in terms - noted $FV(t)$ - and free and bound variables in formulas - noted $FV(\alpha)$ and $BV(\alpha)$, with $FV(\Gamma) = \{FV(\alpha) \ | \ \alpha \in \Gamma\}$;
    \item substitution functions of individual variables with terms in terms - noted $t[x/u]$, $t[x_1, ..., x_n/u_1, ..., u_n]$ for simultaneous substitution defined as usual - and formulas - noted $\alpha[x/t]$, $\alpha[x_1, ..., x_n/t_1, ..., t_n]$ for simultaneous substitution defined as usual - and substitution functions of formulas with formulas in formulas - noted $\alpha[\beta/\gamma]$;
    \item a term being free for a variable in a formula, and a formula being free for a formula in a formula.
\end{itemize}
\end{defn}

\begin{defn}
Given a first-order logical language $L_1$, an \emph{expansion} of $L_1$ is a first-order logical language $L_2$ such that $\texttt{Al}_{L_1} \subseteq \texttt{Al}_{L_2}$, $\texttt{TERM}_{L_1} \subseteq \texttt{TERM}_{L_2}$ and $\texttt{FORM}_{L_1} \subseteq \texttt{FORM}_{L_2}$.
\end{defn}

\subsection{Bases}

The notion of atomic base requires the further introduction of a prior notion of atomic system. 

\subsubsection{Atomic systems}

Following Prawitz, atomic systems will be understood as Post-systems. Rules will be presented in a Gentzen format.

\begin{defn}
An \emph{atomic system} is a pair $\langle L, \mathfrak{R} \rangle$, where $L$ is a first-order logical language, and $\mathfrak{R}$ is a recursively axiomatisable set of rules of the form

\begin{prooftree}
\AxiomC{$\alpha_1, ..., \alpha_n$}
\UnaryInfC{$\beta$}
\end{prooftree}
such that:

\begin{itemize}
    \item $\alpha_i \in \texttt{ATOM}_L$ and $\alpha_i \neq \bot$ ($i \leq n$);
    \item  $\beta \in \texttt{ATOM}_L$ or $\beta = \bot$ and, if $x \in FV(\beta)$, then there is $i \leq n$ such that $x \in FV(\alpha_i)$.
\end{itemize}
\end{defn}
\noindent We leave unspecified whether atomic rules can or cannot bind assumptions and variables, although this is a crucial point when dealing with completeness issues [Piecha, de Campos Sanz \& Schroeder-Heister 2015, Piecha \& Schroeder-Heister 2018].

We quickly outline two examples of atomic system that, together with the next definition, will make somewhat more precise the idea of an atomic system determining the meaning of the non-logical signs of a first-order language. The first example is the \emph{empty system} $\langle \emptyset, \emptyset \rangle$. Another is an atomic system $\texttt{Ar}$ for first-order arithmetic, where we have standard rules for equality plus proper arithmetic rules:

\begin{prooftree}
\AxiomC{}
\RightLabel{($=_R$)}
\UnaryInfC{$t = t$}
\AxiomC{$t = u$}
\RightLabel{($=_S$)}
\UnaryInfC{$u = t$}
\AxiomC{$t = u$}
\AxiomC{$u = z$}
\RightLabel{($=_T$)}
\BinaryInfC{$t = z$}
\noLine
\TrinaryInfC{}
\end{prooftree}

\begin{prooftree}
\AxiomC{$0 = s(t)$}
\RightLabel{($s_1$)}
\UnaryInfC{$\bot$}
\AxiomC{$s(t) = s(u)$}
\RightLabel{($s_2$)}
\UnaryInfC{$t = u$}
\noLine
\BinaryInfC{}
\end{prooftree}

\begin{prooftree}
\AxiomC{}
\RightLabel{($+_1$)}
\UnaryInfC{$t + 0 = t$}
\AxiomC{}
\RightLabel{($+_2$)}
\UnaryInfC{$t + s(u) = s(t + u)$}
\noLine
\BinaryInfC{}
\end{prooftree}

\begin{prooftree}
\AxiomC{}
\RightLabel{($\cdot_1$)}
\UnaryInfC{$t \cdot 0 = 0$}
\AxiomC{}
\RightLabel{($\cdot_2$)}
\UnaryInfC{$t \cdot s(u) = (t \cdot u) + t$}
\noLine
\BinaryInfC{}
\end{prooftree}

\begin{defn}
An atomic system $\texttt{S} = \langle L_1, \mathfrak{R} \rangle$ is \emph{an atomic system for} $L_2$ iff $L_2$ is an expansion of $L_1$.
\end{defn}
\noindent Whatever $L$ is, the empty atomic system is an atomic system for $L$ stating that no deductive relationship holds between the logical signs of $L$. $\texttt{Ar}$ instead exemplifies the idea that an atomic system determines the ‘‘arithmetic meaning" of the non-logical signs of its language for arithmetic.

\begin{defn}
Given an atomic system $\texttt{S}_1 = \langle L_1, \mathfrak{R}_1 \rangle$, we say that an \emph{expansion} of $\texttt{S}_1$ is an atomic system $\texttt{S}_2 = \langle L_2, \mathfrak{R}_2 \rangle$ with $L_2 = L_1$ and $\Re_1 \subseteq \Re_2$, or with $L_2$ proper expansion of $L_1$ and $\Re_1 \subset \Re_2$.
\end{defn}
\noindent Observe that the notion of expansion $\texttt{S}_2$ of an atomic system $\texttt{S}_1$ is here given in a standard set-theoretic way, i.e. $\texttt{S}_2$ contains at least all the rules of $\texttt{S}_1$. If one assumes that the inductive definition of the derivations-set of atomic systems remains ‘‘stable" throughout expansions, this also implies that the derivations-set of $\texttt{S}_2$ includes at least the derivations-set of $\texttt{S}_1$. Otherwise, one can explicitly require this fact, and even adopt a weaker notion of expansion of an atomic system, i.e. $\Gamma \vdash_{\texttt{S}_1} \alpha \Rightarrow \Gamma \vdash_{\texttt{S}_2} \alpha$. In what follows we shall leave aside the discussion of these different formulations, as they are not crucial with respect to the framework we are going to propose. We remark, however, that the choice of one of them becomes central when investigating other issues - e.g., completeness [see Piecha 2016].

\subsubsection{Atomic bases}

\begin{defn}
Given a first-order logical language $L$, let us indicate: with $\texttt{C}$ the set of its individual constants; with $\texttt{R}$ the set of its relational constants; with $\texttt{S}$ an atomic system for $L$. We say that $\langle \texttt{C}, \texttt{R}, \texttt{S} \rangle$ is an \emph{atomic base} on $L$. We call $L$ the \emph{background language} of the base. If $\texttt{S} = \emptyset$, we say that the atomic base is \emph{logical}.
\end{defn}

\begin{defn}
Let $\mathfrak{B}$ be an atomic base on $L_1$ with atomic system $\texttt{S}_1$. We say that an \emph{expansion} of $\mathfrak{B}$ is an atomic base on $L_2$ with atomic system $\texttt{S}_2$, for $L_2$ expansion of $L_1$ and $\texttt{S}_2$ expansion of $\texttt{S}_1$.
\end{defn}

\subsection{Operational types}

We finally introduce the notion of operational type over a background language. The notion of type is a particular case of that of operational type, that obtains when the domain is empty and the co-domain consists of a closed formula.

\begin{defn}
Let $L$ be a first-order logical language. A \emph{pre-type} on $L$ is any $\alpha \in \texttt{FORM}_L$, or an expression
\begin{center}
    $\Gamma \rhd \alpha$
\end{center}
where $\Gamma$ is a finite non-empty sub-set of $\texttt{FORM}_L$ and $\alpha \in \texttt{FORM}_L$. We call $\Gamma$ the \emph{domain} of the pre-type, $\alpha$ the \emph{co-domain} of the pre-type. Given a sequence of pre-types, its \emph{domain} is the union set of the domains of its elements. An \emph{operational type} on $L$ is either a pre-type on $L$ or an expression of the form

\begin{center}
$\tau_1, ..., \tau_n \rhd \tau_{n + 1}$
\end{center}
where $\tau_i$ ($i \leq n + 1$) is a pre-type and the domain of $\tau_{n + 1}$ is a subset of the domain of $\tau_1, ..., \tau_n$ - in the case when the operational type is a closed $\alpha \in \texttt{FORM}_L$, we shall simply speak of a \emph{type} on $L$. We say that $\tau_1, ..., \tau_n$ is the \emph{domain} of the operational type, and that $\tau_{n + 1}$ is the \emph{co-domain} of the operational type. We say that each $\tau_i$ ($i \leq n$) is an \emph{entry} of the domain operational type.
\end{defn}
\noindent As a notational convention, if the domain of an operational type is non-empty, and some entry of the domain or the co-domain have non-empty domain, we put domain and co-domain between parentheses.\footnote{Observe that expressions like the following are not operational types: $(\rhd \alpha) \rhd \beta$, $(\rhd \alpha) \rhd (\rhd \beta)$, $(\Gamma \rhd \alpha) \rhd (\rhd \beta)$, $\rhd \beta$. This depends on the requirement that $\Gamma$ be a \emph{non-empty} set of formulas, so that $\rhd \alpha$ and $\rhd \beta$ are not pre-types. Our choice is motivated simply by the conceptual order we decided to adopt for discussing operations on grounds as shown in Section 4.2. This order classifies operations on grounds, first, depending on whether their domain is or not empty and, second, according to the ‘‘complexity" of the elements of the domain - i.e. according to whether all the elements of the domain are formulas, or there is an element of the form $\Gamma \rhd \alpha$. An alternative, but equivalent account can of course be given, in such a way that $\rhd \alpha$ and $\rhd \beta$ become pre-types.}

\section{Grounds and operations on grounds}

Here, we define general notions of grounds and operations on grounds. The languages for speaking of these objects are introduced in Section 5, and the denotation functions through which this ‘‘speaking of" obtains are introduced in Section 6. The notions of ground and operation on grounds run by simultaneous recursion, based on the concept of ground for categorical assertions.

\subsection{Grounds}

Grounds for categorical assertions are given through clauses where \emph{primitive} operations are assumed, one for each first-order logical constant. Primitive operations stand in need of no definition, as the clauses make them meaning-constitutive. Moreover - as any other operation of our theory - they come with an operational type. However, we leave operational types implicit, except when it is necessary to indicate them - and in these cases, we indicate them in square brackets after the symbol for the operation. Some primitive operations bind individual or typed-variables, which we indicate by accompanying the symbol for the operation with the variable that is bound. All the formulas in the clauses have to be understood as closed, or as open but then closed with adequate bindings. So, given an atomic base $\mathfrak{B}$ on a first-order language $L$ with atomic system $\texttt{S}$ we say that:
\begin{itemize}
    \item[($\texttt{At}$)] a ground on $\mathfrak{B}$ for $\vdash \alpha$ with $\alpha \in \texttt{ATOM}_L$ is any derivation of $\alpha$ in $\texttt{S}$ with no undischarged assumption or unbound individual variables;
    \item[($\wedge$)] $g_1$ is a ground on $\mathfrak{B}$ for $\vdash \alpha$ and $g_2$ is a ground on $\mathfrak{B}$ for $\vdash \beta$ iff $\wedge I(g_1, g_2)$ is a ground on $\mathfrak{B}$ for $\vdash \alpha \wedge \beta$;
    \item[($\vee$)] $g$ is a ground on $\mathfrak{B}$ for $\vdash \alpha_i$ iff $\vee I [\alpha_i \rhd \alpha_1 \vee \alpha_2] (g)$ is a ground on $\mathfrak{B}$ for $\vdash \alpha_1 \vee \alpha_2$ ($i = 1, 2$);
    \item[($\rightarrow$)] $f(\xi^\alpha)$ is a ground on $\mathfrak{B}$ for $\alpha \vdash \beta$ iff $\rightarrow I \xi^\alpha (f(\xi^\alpha))$ is a ground on $\mathfrak{B}$ for $\vdash \alpha \rightarrow \beta$;
    \item[($\exists$)] $g$ is a ground on $\mathfrak{B}$ for $\vdash \alpha(t)$ iff $\exists I [\alpha(t) \rhd \exists x \alpha(x)](g)$ is a ground on $\mathfrak{B}$ for $\vdash \exists x \alpha(x)$;
    \item[($\forall$)] $f(x)$ is a ground on $\mathfrak{B}$ for $\alpha(x)$ iff $\forall I x (f(x))$ is a ground on $\mathfrak{B}$ for $\vdash \forall x \alpha(x)$.
\end{itemize}
There is finally a clause like ‘‘nothing else is a closed ground on $\mathfrak{B}$", saying that these are the only possible forms of closed grounds on $\mathfrak{B}$. As for $\bot$, we put:
\begin{itemize}
    \item[($\bot$)] there is no ground on $\mathfrak{B}$ for $\vdash \bot$.
\end{itemize}
On every atomic base $\mathfrak{B}$ on $L$ we can thus introduce a $\mathfrak{B}$-operation on grounds of operational type $\bot \rhd \alpha$ for every $\alpha \in \texttt{FORM}_L$ defined by an empty set of computation instructions, i.e. the empty function of type $\bot \rhd \alpha$. The clauses ($\rightarrow$) and ($\forall$) involve the so far undefined notions of, respectively, ground for $\alpha \vdash \beta$ and ground for $\alpha(x)$. These notions require the introduction of the more general concept of operation on grounds.

\subsection{First- and second-level operations}

Operations on grounds will be generally indicated by a notation of the form
\begin{center}
$f(x_1, ..., x_n, \xi^{\tau_1}, ..., \xi^{\tau_m})$
\end{center}
for $n, m \geq 0$, where $\tau_i$ is a possibly open formula, or a hypothetical or general-hypothetical assertion ($i \leq m$). So, the operation at issue is defined on $n$ individuals $k_1, ..., k_n$, and on grounds for $\tau_i[k_1, ..., k_n/x_1, ..., x_n]$ - plus possibly other replacements with arbitrary individuals of unbound variables not occurring in $\tau_i$ other than $x_1, ..., x_n$ ($i \leq m$).

Let $\mathfrak{B}$ be an atomic base on a first-order language $L$. We assume that $\mathfrak{B}$ is associated to a domain $D_\mathfrak{B}$, for example by specifying (in a possibly inductive way) some canonical form that the elements of $D_\mathfrak{B}$ must have - e.g., for first-order arithmetic, we could give formation rules in a Martin-L\"{o}f style $0 \in \mathbb{N}$ and $n \in \mathbb{N} \Rightarrow s(n) \in \mathbb{N}$. We say that a $\mathfrak{B}$-operation on grounds of operational type
\begin{center}
$\alpha(x_1, ..., x_n)$
\end{center}
with $n \geq 0$, is a total constructive function
\begin{center}
$f(x_1, ..., x_m)$
\end{center}
with $m \geq n$ and $m > 0$ such that, for every $k_1, ..., k_m \in D_\mathfrak{B}$,
\begin{center}
$f(k_1, ..., k_m)$
\end{center}
is a ground for
\begin{center}
$\vdash \alpha(k_1, ..., k_m/x_1, ..., x_n)$
\end{center}
- what is replaced in the formula must be a name of individual, but we will not indicate this neither here nor in what follows. The operation is a ground on $\mathfrak{B}$ for the corresponding general assertion, proper if $m = n$, improper otherwise. As an example, think of the following reasoning:
\begin{enumerate}
    \item $x > y$, assumption
    \item $y > z$, assumption
    \item hence $x > z$, from 1 and 2, $>$ being transitive
    \item hence, $y > z \rightarrow x > z$, from 3, introducing implication on 1
    \item hence, $x > y \rightarrow (y > z \rightarrow x > z)$, from 4, introducing implication on 2
    \item $w = w \rightarrow (x > y \rightarrow (y > z \rightarrow x > z))$, from 5, introducing implication
    \item $w = w$, $=$ being reflexive
    \item hence, $x > y \rightarrow (y > z \rightarrow x > z)$, from 6 and 7, by \emph{modus ponens}
\end{enumerate}
If in this reasoning we replace $x, y, z, w$ with four natural numbers $n, m, p, q$, we obtain a correct reasoning proving $n > m \rightarrow (m > p \rightarrow n > p)$. Observe that $w$ only occurs in steps 6 and 7, which are patently redundant. If we eliminate them and make step 8 correspond to step 5, the reduced reasoning thereby obtained still proves our conclusion. This redundancy explains why the expression of an operation on grounds may involve more individual variables than those that occur in its type - accordingly, in the reduced reasoning there are no more individual variables than those that occur in the conclusion. Similarly, we shall say that a $\mathfrak{B}$-operation on grounds of operational type

\begin{center}
    $\alpha_1(\underline{x}_1), ..., \alpha_n(\underline{x}_n) \rhd \beta(\underline{x}_{n + 1})$
\end{center}
is a total constructive function
\begin{center}
    $f(\underline{x}_1, ..., \underline{x}_{n + 1}, \xi^{\alpha_1(\underline{x}_1)}, ..., \xi^{\alpha_n(\underline{x}_n)})$
\end{center}
such that, for every $\underline{k}$, for every $g_i$ ground on $\mathfrak{B}$ for $\vdash \alpha_i(\underline{k}/\underline{x}_i)$ ($i \leq n$),

\begin{center}
    $f(\underline{k}, g_1, ..., g_n)$
\end{center}
is a ground on $\mathfrak{B}$ for $\vdash \beta(\underline{k}/\underline{x}_{n + 1})$. Operations of this kind may also bind individual variables on some entries of index $i$ - although these individual variables may occur free in entries of index $j \neq i$ without being bound on these entries. If an individual variable is bound on all the entries where it occurs, then we should consider it as no longer belonging to the sequences $\underline{x}_1, ..., \underline{x}_{n + 1}$ in the notation that indicates the operation. Also, we should take into account the possibility of having a sequence $\underline{z}$ of vacuously bound individual variables, as well as the possibility of having a sequence $\underline{y}$ of individual variables not occurring in the operational type. For $\underline{y}$, the following restriction holds: if, for every $i \leq n$, the operation binds some $x$ on entries of index $i$, then $\underline{y}$ is empty. The reason for this restriction is seen from the example with ($\forall I$) below. The operation is a ground for the corresponding hypothetical or general-hypothetical assertion, proper if $\underline{y}$ is empty, improper otherwise. As for the examples, consider the following reasoning:

\begin{enumerate}
    \item $2 + 2 = 4$, assumption
    \item $4 = x$, assumption
    \item hence $2 + 2 = x$, from $1$ and $2$, $=$ being transitive
\end{enumerate}
If in this reasoning we replace $x$ with $\sqrt{16}$, and if we replace the assumptions with computations of $2 + 2 = 4$ and of $\sqrt{16} = 4$, we obtain a correct reasoning proving $2 + 2 = \sqrt{16}$. Observe that no variable is bound in this case. Now consider $\forall$-introduction in two variables in Gentzen's first-order natural deduction
\begin{prooftree}
\AxiomC{$\alpha(x, y)$}
\RightLabel{($\forall_I$)}
\UnaryInfC{$\forall x \alpha(x, y)$}
\end{prooftree}
Clearly, we cannot consider this as an \emph{argument}, for otherwise $\alpha(x, y)$ would be an assumption, and it would violate the restriction on the so-called proper parameter. But still we can take it - and one usually takes it - as a \emph{reasoning rule}, and hence as an operation $\forall I x(y, \xi^{\alpha(x, y)})$ that, binding $x$, yields grounds for $\vdash \forall x \alpha(x, k)$ from grounds for $\vdash \alpha(x, k/y)$, i.e. from operations on grounds of operational type $\alpha(x, k/y)$. Since $\alpha(x, y)$ cannot be considered as an assumption, no individual variables other than $y$ can be taken into account. This explains the restriction on $\underline{y}$ above.

We can now describe $\mathfrak{B}$-operations on grounds of second-level. Recall that in this case the entries of the operational type may be general-hypothetical assertions of first-level. With $\underline{x}^\downarrow$ we indicate the sequence of the individual variables occurring free in an entry, if this entry is a simple formula, or in the \emph{co-domain} of this entry, if it is a general-hypothetical assertion of first-level. A $\mathfrak{B}$-operation on grounds of operational type
\begin{center}
    $\tau_1, ..., \tau_n \rhd \tau_{n + 1}$
\end{center}
is a total constructive function
\begin{center}
    $f(\underline{x}^\downarrow_1, ..., \underline{x}^\downarrow_{n + 1}, \xi^{\tau_1}, ..., \xi^{\tau_n})$
\end{center}
such that, for every $\underline{k}$, for every $g_i$ ground on $\mathfrak{B}$ for $\vdash \tau_i[\underline{k}/\underline{x}_i]$,
\begin{center}
    $f(\underline{k}, g_1, ..., g_n)$
\end{center}
is a ground on $\mathfrak{B}$ for $\vdash \tau_{n + 1}[\underline{k}/\underline{x}_{n + 1}]$.\footnote{We should also allow the operation to produce a ground for the co-domain of the co-domain when applied simply to grounds for the co-domains of the domain, i.e. called $\vdash \alpha_i(\underline{x}_i)$ the co-domain of the $i$-th entry and $\beta(\underline{x}_{n + 1})$ the co-domain of the co-domain, for every $\underline{k}$, for every $g_i$ ground on $\mathfrak{B}$ for $\vdash \alpha_i(\underline{x}_i)[\underline{k}/\underline{x}_i]$, $f(\underline{k}, g_1, ..., g_n)$ is a ground on $\mathfrak{B}$ for $\vdash \beta(\underline{k}/\underline{x}_{n + 1})$. We have omitted this case for greater simplicity, but it becomes central when proving the denotation theorem of Section 6.2.} These operations may again bind individual variables, and require sequences $\underline{z}$ for vacuously bound individual variables or $\underline{y}$ for individual variables not occurring in the operational type. We have in this case the same wording as for operations of first-level. Also, these operations may bind typed-variables on some entries of index $i$, corresponding to dischargement of the formula constituting the type of the variable on the corresponding entry - although these formulas may remain undischarged in entries of index $j \neq i$. Finally, we may need to add a sequence of vacuously bound typed-variables. We remark that typed-variables are bound under replacement of individual variables with $\underline{k}$ in the types of the original bound typed-variables. The operation is to be a ground on $\mathfrak{B}$ for the corresponding hypothetical or general-hypothetical assertion, proper for $\underline{y}$ empty, improper otherwise. Observe that the operation is defined on individuals that correspond only to individual variables in the co-domains of the entries. The reasons for this fact will be soon clear. A typical example of operation on grounds of second-level is $\rightarrow$-introduction in Gentzen's natural deduction, that we consider here in the following form:
\begin{prooftree}
\AxiomC{$[\alpha(x)]$}
\noLine
\UnaryInfC{$\vdots$}
\noLine
\UnaryInfC{$\beta(y)$}
\RightLabel{($\rightarrow_I$)}
\UnaryInfC{$\alpha(x) \rightarrow \beta(y)$}
\end{prooftree}
This is an operation on grounds $\rightarrow I \xi^{\alpha(x)} (x, y, \xi^{\alpha(x) \rhd \beta(y)})$ of operational type $(\alpha(x) \rhd \beta(y)) \rhd \alpha(x) \rightarrow \beta(y)$, binding $\xi^{\alpha(x)}$: for every $t_1, t_2$, for every ground $g$ for $\alpha(t_1) \vdash \beta(t_2)$, or simply for $\vdash \beta(t_2)$, $\rightarrow I \xi^{\alpha(t_1)} (t_1, t_2, g)$ is a ground for $\vdash \alpha(t_1) \rightarrow \beta(t_2)$.

$x$ is mentioned \emph{only because} $\alpha(x)$ \emph{occurs in the co-domain}: as we shall see in Section 4.4, the operation can be applied to operations on grounds $h(x, y, \xi^{\alpha(x)})$ of operational type $\alpha(x) \rhd \beta(y)$, so that the application will anyway be defined on individuals corresponding to $x$, i.e. $f(y, h(x, y, \xi^{\alpha(x)}))$. Moreover, as we shall see in Section 4.4, the operation may be applied simply to grounds for $\vdash \beta(t)$, for every $t$, or to operations on grounds $h(y)$ of operational type $\beta(y)$, or finally to operations on grounds whose operational type involves \emph{neither} $\alpha(x)$ \emph{nor} $x$. So, $x$ should not be mentioned as referred to the occurrence of $\alpha(x)$ in the domain of the domain.\footnote{In order not to burden the notation excessively, some details have been omitted. For example, we have not mentioned the case of operations on grounds of second-level $f$ that \emph{vacuously} bind $\xi^\alpha$, for $\alpha$ occurring in the domain of some entry of their operational type. This would have required to consider the typed-variables of the linguistic expression of an operation on grounds of first-level $h$ as coming with an index, and $f$ to bind typed-variables \emph{of a specific index} $i$. The binding is then vacuous when $\xi^\alpha$ has in $h$ an index $j \neq i$. Also, we have not defined what we may call \emph{variables-variants}. The variable-variant of an operation on grounds of operational type $\alpha_1(\underline{x}_1), ..., \alpha_m(\underline{x}_n) \rhd \beta(\underline{x}_{n + 1})$ is to be an operation on grounds of the same operational type, except that it may take as values, not grounds for the entries of the operational type after replacement of all the free individual variables with (names of) individuals, but operation on grounds for these very entries where some individual variables may remain free. We will not really need these notions in what follows, so we have dealt with them in this footnote for the sake of completeness.}

Finally, observe that the notion of $\mathfrak{B}$-operation on grounds as defined here is non-monotonic, i.e. something may be a $\mathfrak{B}$-operation on grounds without being a $\mathfrak{B}^+$-operation on grounds for some expansion $\mathfrak{B}^+$ of $\mathfrak{B}$. To see this, let $\mathfrak{B}$ be a logical base on a standard background language for propositional logic. The empty function of operational type $p \rhd \bot$ is a $\mathfrak{B}$-operation on grounds of operational type $p \rhd \bot$, simply because there are no grounds for $p$ on $\mathfrak{B}$. Take now the expansion $\mathfrak{B}^+$ of $\mathfrak{B}$ on the same background language as $\mathfrak{B}$, and whose atomic system is $\{p\}$ - i.e. has $p$ as unique axiom. Our empty function is no longer a $\mathfrak{B}^+$-operation on grounds, for we now have a ground for $p$, but $\mathfrak{B}^+$ is consistent. More generally, there is \emph{no} $\mathfrak{B}^+$-operation on grounds of operational type $p \rhd \bot$. This also shows that being a ground on a logical base does not amount to being a ground on every base [if instead one wants this to hold, one must require that $\mathfrak{B}$-operations on grounds are ‘‘ground-preserving", not only over $\mathfrak{B}$, but also over expansions of $\mathfrak{B}$, i.e. that they produce $\mathfrak{B}^+$-grounds for the co-domain from $\mathfrak{B}^+$-grounds for the domain, for every expansion $\mathfrak{B}^+$ of $\mathfrak{B}$ - one would thereby obtain a ground-theoretic version of what Schroeder-Heister 2006 does for Prawitz's semantics of valid arguments].

\subsection{Defining equations}

Since we expect operations on grounds to be \emph{effective} functions, we must require them to be \emph{defined} by one or more equations, showing \emph{how} they produce given results on given arguments.

We shall impose no restriction on the form of a defining equation, except requiring that the instructions it provides must render total and effective the operation it is intended to define. For the examples, we here limit ourselves to closed formulas, and introduce equations for operations on grounds corresponding to Gentzen's eliminations in first-order natural deduction, as well as to induction in first-order arithmetic:
\begin{itemize}
    \item $f_\wedge(\xi^{\alpha_1 \wedge \alpha_2})$ of operational type $\alpha_1 \wedge \alpha_2 \rhd \alpha_i$ defined by the equation
    \begin{center}
        $f_\wedge(\wedge I (g_1, g_2)) = g_i$
    \end{center}
    ($i = 1, 2$);
    \item $f_\vee(\xi^{\alpha_1 \vee \alpha_2}, \xi^{\alpha_1 \rhd \beta}, \xi^{\alpha_2 \rhd \beta})$ of operational type $\alpha_1 \vee \alpha_2, (\alpha_1 \rhd \beta), (\alpha_2 \rhd \beta) \rhd \beta$ binding $\xi^{\alpha_i}$ on the $i + 1$th entry defined by the equation
    \begin{center}
        $f_\vee(\vee I[\alpha_i \rhd \alpha_1 \vee \alpha_2](g), h_1(\xi^{\alpha_1}), h_2(\xi^{\alpha_2})) = h_i(g)$
    \end{center}
    ($i = 1, 2$);
    \item $f_\rightarrow(\xi^{\alpha \rightarrow \beta}, \xi^{\alpha})$ of operational type $\alpha \rightarrow \beta, \alpha \rhd \beta$ defined by the equation
    \begin{center}
        $f_\rightarrow(\rightarrow I \xi^\alpha (h(\xi^\alpha)), g) = h(g)$;
    \end{center}
    \item $f_\forall(\xi^{\forall x \alpha(x)})$ of operational type $\forall x \alpha(x) \rhd \alpha(k)$ defined by the equation
    \begin{center}
        $f_\forall(\forall I x (h(x)) = h(k)$;
    \end{center}
    \item $f_\exists(\xi^{\exists x \alpha(x)}, \xi^{\alpha(x) \rhd \beta})$ of operational type $\exists x \alpha(x), (\alpha(x) \rhd \beta) \rhd \beta$ binding $x$ and $\xi^{\alpha(x)}$ on the second entry defined by the equation
    \begin{center}
        $f_\exists(\exists I[\alpha(t) \rhd \exists x \alpha(x)](g), h(x, \xi^{\alpha(x)})) = h(t, g)$;
    \end{center}
    \item $f^n_{\texttt{Ind}}(\xi^{\alpha(0)}, \xi^{\alpha(x) \rhd \alpha(s(x))})$ of operational type $\alpha(0), (\alpha(x) \rhd \alpha(s(x))) \rhd \alpha(n)$ binding $x$ and $\xi^{\alpha(x)}$ on the second entry defined by the equation
    \begin{center}
        $f^n_{\texttt{Ind}}(g, h(x, \xi^{\alpha(x)})) = \begin{cases} g & \text{if} \ n = 0 \\ h(n - 1, f^{n - 1}_{\texttt{Ind}}(g, h(x, \xi^{\alpha(x)})) & \text{if} \ n > 0 \end{cases}$
    \end{center}
\end{itemize}
In the open case, we specify the behavior of the operation after replacement of individual variables with (names of) individuals and application to grounds for formulas closed by replacement. It is more important to observe that the symbol $=$ \emph{does not} indicate a relation between objects in our semantic universe, for $f_\wedge(\wedge I(g_1, g_2))$ \emph{is not} an object in this universe. The symbol $=$ is, so to say, on a meta-level, and simply shows the result we obtain by applying functions to arguments.

\subsection{Composition and restrictions}

We allow that an operation on grounds can be plugged into another operation on grounds provided the operational type of the former has (co-domain of the) co-domain identical to (the co-domain of) one of the entries of the operational type of the latter.

More specifically, given a base $\mathfrak{B}$, given $h_i(\underline{x}_i, \xi^{\tau^i_1}, ..., \xi^{\tau^i_{m_i}})$ ground on $\mathfrak{B}$ for $\tau^1_i, ..., \tau^i_{m_i} \vdash \tau^i_{m_i + 1}$, where $\tau^i_{m_i + 1}$ is $\alpha_i$ if $\tau^i_{m + 1}$ has empty domain, or it has co-domain $\alpha_i$ otherwise, and given a $\mathfrak{B}$-operation on grounds $f(\underline{x}, \xi^{\tau_1}, ..., \xi^{\tau_n})$ of operational type $\tau_1, ..., \tau_n \rhd \tau_{n + 1}$, where $\tau_i$ is $\alpha_i$ if $\tau_i$ has empty domain, or it has co-domain $\alpha_i$ otherwise, we say that
\begin{center}
    $f(\underline{x}, \dots, h_i(\underline{x}_i, \xi^{\tau^i_1}, ..., \xi^{\tau^i_{m_i}}), \dots)$
\end{center}
is a composition of
\begin{center}
    $f(\underline{x}, \xi^{\tau_1}, ..., \xi^{\tau_n})$ with $h_i(\underline{x}_i, \xi^{\tau^i_1}, ..., \xi^{\tau^i_{m_i}})$
\end{center}
- for the sake of simplicity, here and in the sequel we simply speak of composition of $f$ with $h_i$ - and also that $h_i$ is plugged into $f$ on index $i$ ($i \leq n$). Of course, if $f$ binds $x$ on index $i$, and if $x$ is an element of $\underline{x}_i$, $x$ is to be taken as bound on $i$ in the composition of $f$ with $h_i$; similarly, if $f$ binds $\xi^\alpha$ on index $i$, then every $\xi^{\tau^i_j} = \xi^\alpha$ ($j \leq m_i$) is to be taken as bound on index $i$ in the composition of $f$ with $h_i$.

It is in general not immediate that the composition of grounds yields grounds. To see this, suppose we are on a base $\mathfrak{B}$ for first-order arithmetic as in Section 3.2.1, i.e. with atomic system $\texttt{Ar}$ - and whose domain we indicate as usual with $\mathbb{N}$. Let us consider a primitive $\mathfrak{B}$-operation on grounds $\forall I x(\xi^{\alpha(x)})$ of operational type $x = 0 \rhd \forall x (x = 0)$, and the identity $\mathfrak{B}$-operation on grounds $\texttt{Id}(x, \xi^{x = 0})$ of operational type $x = 0 \rhd x = 0$. Let us now consider the composition of these two operations, i.e. $\forall I x (\texttt{Id}(x, \xi^{x = 0}))$. This should be a $\mathfrak{B}$-operation on grounds again of operational type $x = 0 \rhd \forall x (x = 0)$, but in fact it is not. If it were, then for every $k \in \mathbb{N}$, for any closed derivation $\Delta$ in $\texttt{Ar}$ for $k = 0$, $\forall I x (\Delta)$ should be a ground on $\mathfrak{B}$ for $\forall x (x = 0)$ and, by clause ($\forall$) in Section 4.1, this can be only if $\Delta$ is in turn a $\mathfrak{B}$-operation on grounds of operational type $x = 0$, i.e. if for every $k \in \mathbb{N}$, $\Delta[k/x]$ is an atomic derivation in $\texttt{Ar}$ for $k = x$. If we let $k$ be $0$, and $\Delta$ be an atomic derivation in $\texttt{Ar}$ for $0 = 0$, then for every $k$ other than $0$, say $1$, since $\Delta$ is closed we have $\Delta[1/x] = \Delta$, so $\Delta[1/x]$ cannot be an atomic derivation in $\texttt{Ar}$ for $1 = 0$, so $\forall I x (\Delta)$ is not a ground on $\mathfrak{B}$ for $\forall I x (x = 0)$, so $\forall I x (\texttt{Id}(x, \xi^{x = 0}))$ is not a $\mathfrak{B}$-operation on grounds of operational type $x = 0 \rhd \forall x (x = 0)$.
 
Therefore, we must introduce two restrictions - which are nothing but a generalization of those for quantifiers in Gentzen's natural deduction. For every base $\mathfrak{B}$, for every $\mathfrak{B}$-operation on grounds $f$ of operational type $\tau_1, ..., \tau_n \rhd \tau_{n + 1}$:
\medskip

\noindent \textbf{restriction 1} - if, for some $i \leq n$, $\tau_i$ has empty domain and is of the form
\begin{center}
    $\dots \forall x \dots \alpha(x, t/y)$
\end{center}
or it has non-empty domain, and its co-domain has such form, and if $\tau_{n + 1}$ has empty domain and is of the form
\begin{center}
    $\dots \exists y \dots \beta(u/x, y)$
\end{center}
or it has non-empty domain, and its co-domain has such form, then $t$ is free for $y$ in $\dots \forall x \dots \alpha(x, y)$ and $u$ is free for $x$ in $\dots \exists y \dots \beta(x, y)$.
\medskip

\noindent \textbf{restriction 2} - if $f$ binds $x$ on index $i$, then

\begin{itemize}
    \item every ground $h_i$ on $\mathfrak{B}$ for $\tau^i_1, ... \tau^i_{m + 1} \vdash \tau^i_{m + 1}$ plugged into $f$ on index $i$ must be such that, for every $j \leq m_i$, if $\tau^i_j$ has non-empty domain $\Gamma$ then, for every $\beta \in \Gamma$, $x$ occurs free in $\beta$ iff $h_i$ binds $\xi^\beta$ on index $^i_j$, or $h_i$ is in turn a composed operation and $x$ is bound on index $^i_j$ by some operation which $h_i$ is composed of, or $f$ binds $\xi^\beta$ on index $i$ - of course, if $\tau^i_j$ ha empty domain, then it may be taken to be simply $\beta$ - and
    \item $x$ does not occur free in $\tau_{n + 1}$, if $\tau_{n + 1}$ has empty domain, or in the co-domain of $\tau_{n + 1}$, if $\tau_{n + 1}$ has non-empty domain.
\end{itemize}

\begin{prp}
Let $\mathfrak{B}$ be an atomic base, and let $f(\underline{x}, \xi^{\tau_1}, ..., \xi^{\tau_n})$ be a $\mathfrak{B}$-operation on grounds of operational type $\tau_1, ..., \tau_n \rhd \tau_{n + 1}$. Then, for every $h_i(\underline{x}_i, \xi^{\tau^i_1}, ..., \xi^{\tau^i_{m_i}})$ ground on $\mathfrak{B}$ for $\tau^i_1, ... \tau^i_{m + 1} \vdash \tau^i_{m + 1}$ ($i \leq n$), the composition of $f$ with $h_i$ is a ground on $\mathfrak{B}$ for
\begin{center}
    $\dots \{\tau^i_1, ..., \tau^i_{m_i}\} - \nu_i \dots \vdash (\dots ((\{\sigma^i_s, ..., \sigma^i_t\} - \{\nu^i_s, ..., \nu^i_t\}) - \nu_i) \dots \vdash \beta)$
\end{center}
where
\begin{itemize}
    \item $\nu_i$ is the set of the types of the typed-variables bound by $f$ on index $i$;
    \item $\{\nu^i_s, ..., \nu^i_t\}$ is the set of the types of the typed-variables bound by $h_i$ on indexes $^i_s, ..., ^i_t$ ($s, ..., t \leq m_i$)
    \item $\{\sigma^i_s, ..., \sigma^i_t\}$ is the the set of the domains of the entries of the operational type of $h_i$ with non-empty domain ($s, ..., t \leq m_i$);
    \item $\beta$ is $\tau_{n + 1}$, if $\tau_{n + 1}$ has empty domain, or it is the co-domain of $\tau_{n + 1}$, if $\tau_{n + 1}$ has non-empty domain.
\end{itemize}
Moreover, if the composition of $f$ with $h_i$ is defined on $v \geq 0$ individuals, for any sequence $\underline{k}$ of individuals of length $w \leq v$, the composition of $f$ with $h_i$ applied to $\underline{k}$ is a ground on $\mathfrak{B}$ for the same assertion as above, with appropriate replacements.
\end{prp}

\subsection{Equivalence and identity}

We finally introduce equivalence and identity between grounds on an atomic base $\mathfrak{B}$ with atomic system $\texttt{S}$, to which a domain of individuals $D_\mathfrak{B}$ is associated. Equivalence on $\mathfrak{B}$ is defined by simultaneous recursion, based on the ‘‘categorical" case. Given two grounds $g_1, g_2$ on $\mathfrak{B}$, $g_1 \equiv_\mathfrak{B} g_2$ iff

\begin{itemize}
    \item they have the same operational type and are defined on the same number of individuals of $D_\mathfrak{B}$;
    \item the domain of their operational type is empty and they are defined on a null number of individuals $\Rightarrow$
    \begin{itemize}
        \item $g_1$ and $g_2$ are the same atomic derivation in $\texttt{S}$;
        \item $g_1$ is $\wedge I(g_3, g_4)$ and $g_2$ is $\wedge I(g_5, g_6)$ with $g_3 \equiv_\mathfrak{B} g_5$ and $g_4 \equiv_\mathfrak{B} g_6$;
        \item $g_1$ is $\vee I [\alpha_i \rhd \alpha_1 \vee \alpha_2](g_3)$ and $g_2$ is $\vee I[\alpha_i \rhd \alpha_1 \vee \alpha_2](g_4)$ with $g_3 \equiv_\mathfrak{B} g_4$  ($i = 1, 2$);
        \item $g_1$ is $\rightarrow I \xi^\alpha (g_3)$ and $g_2$ is $\rightarrow I \xi^\alpha (g_4)$ with $g_3 \equiv_\mathfrak{B} g_4$;
        \item $g_1$ is $\exists I [\alpha(t) \rhd \exists x \alpha(x)](g_1)$ and $g_3$ is $\exists I [\alpha(t) \rhd \exists x \alpha(x)](g_4)$ with $g_3 \equiv_\mathfrak{B} g_4$;
        \item $g_1$ is $\forall I x (g_3)$ and $g_2$ is $\forall I x (g_4)$ with $g_3 \equiv_\mathfrak{B} g_4$.
    \end{itemize}
    \item the domain of their operational type is non-empty or they are defined on a non-null number of individuals $\Rightarrow$ when applied to identical individuals of $D_\mathfrak{B}$, and then to identical grounds on $\mathfrak{B}$ for elements in their domain with appropriate substitutions of individual variables with (names of) individuals  (see below for identity), in such a way as to obtain respective grounds $g_1^*$ and $g_2^*$ on $\mathfrak{B}$ whose operational type has empty domain and which are defined on a null number of individuals, it holds that $g_1^*$ and $g_2^*$ are identical (see below for identity).
\end{itemize}
Observe that equivalence is a kind of extensional identity, i.e. two operations on grounds are equivalent when they produce identical values from identical arguments. This characterisation requires of course to say what identity is, and this we do immediately below. Before that, we just remark that we take it to be part of their being \emph{operations} that the following holds of operations on grounds.

\paragraph{Fact} Any $\mathfrak{B}$-operation on grounds $f(x_1, ..., x_n, \xi^{\tau_1(\underline{x}_1)}, ..., \xi^{\tau_m(\underline{x}_m)})$ is defined in such a way that, for every $k_1, ..., k_n \in D_\mathfrak{B}$, for every $g^1_i, g^2_i$ grounds on $\mathfrak{B}$ for elements in the domain of the operational type with appropriate substitutions of individual variables with (names of) individuals, and such that $g^1_i \equiv_\mathfrak{B} g^2_i$ ($i \leq m$),

\begin{center}
    $f(k_1, ..., k_n, g^1_1, ..., g^1_m) \equiv_\mathfrak{B} f(k_1, ..., k_n, g^2_1, ..., g^2_m)$
\end{center}

\noindent In other words, an operation on grounds must be so defined as to yield equivalent values from equivalent arguments.

Identity on $\mathfrak{B}$, indicated with $\approx_\mathfrak{B}$, is defined by distinguishing two cases, depending on whether the objects are atomic derivations of $\texttt{S}$ or made of (possibly composite) operations on grounds. It involves a notion of operations \emph{defined by the same equation}, that we shall take here in the following intuitive sense: since a defining equation provides instructions for computing an operation, two operations are defined by the same equation when their computation instructions are the same. Given two grounds $g_1, g_2$ on $\mathfrak{B}$, $g_1 \approx_\mathfrak{B} g_2$ iff

\begin{itemize}
    \item they have the same operational type and are defined on the same number of individuals of $D_\mathfrak{B}$;
    \item $g_1$ and $g_2$ are the same atomic derivation in $\texttt{S}$;
    \item $g_1$ is $f(x_1, ..., x_n, \sigma_1, ..., \sigma_m)$ and $g_2$ is $h(x_1, ..., x_n, \sigma^*_1, ..., \sigma^*_m)$, where
    \begin{center}
    $f(x_1, ..., x_n, \xi^{\tau_1}, ..., \xi^{\tau_m})$ and $h(x_1, ..., x_n, \xi^{\tau_1 *}, ..., \xi^{\tau_m *})$
    \end{center}
    are defined by the same equation and $\sigma_i$ is a ground on $\mathfrak{B}$ iff $\sigma^*_i$ is a ground on $\mathfrak{B}$, with $\sigma_i \approx_\mathfrak{B} \sigma^*_i$ ($i \leq m$).
\end{itemize}
Identity is therefore ‘‘sameness", and clearly $g_1 \approx_\mathfrak{B} g_2 \Rightarrow g_1 \equiv_\mathfrak{B} g_2$, although the vice versa may not hold. As an example of equivalent but not identical operations on grounds, consider first of all the operation
\begin{center}
    $DS(\xi^{\alpha \vee \beta}, \xi^{\neg \alpha})$
\end{center}
of operational type $\alpha \vee \beta, \neg \alpha \rhd \beta$, defined by requiring that, for every $\vee I (g_1)$ ground for $\vdash \alpha \vee \beta$ with $g_1$ ground for $\vdash \beta$, for every $g_2$ ground for $\vdash \neg \alpha$,
\begin{center}
    $DS(\vee I(g_1), g_2) = g_1$.
\end{center}
It is equivalent, but not identical, to the composite $f_\vee (\xi^{\alpha \vee \beta}, \bot_\beta(f_\rightarrow(\xi^{\alpha}, \xi^{\neg \alpha})), \xi^\beta)$, where $f_\vee$ and $f_\rightarrow$ are the operations defined in Section 4.3, and where $\bot_\beta$ is the empty function of operational type $\bot \rhd \beta$.\footnote{The definition of identity has to be slightly modified if one takes into account vacuous binding of typed-variables occurring in the domain of some entry of an operation on grounds of second-level. For example, $\rightarrow I \xi^{\alpha}_2(\rightarrow I \xi^{\alpha}_1(\xi^{\alpha}_1))$ is not to be identical - nor even equivalent - to $\rightarrow I \xi^{\alpha}_1(\rightarrow I \xi^{\alpha}_2(\xi^{\alpha}_1))$, although both are grounds for $\vdash \alpha \rightarrow (\alpha \rightarrow \alpha)$. This can be done by taking into account what said in footnote 1, and adding that for two operations on grounds of second-level to be identical, there must be a ‘‘structure-preserving" re-indexing of the typed-variables such that all the typed-variables have the same indexes in the two cases - e.g. $\rightarrow I \xi^{\alpha}_2(\rightarrow I \xi^{\alpha}_1(\xi^{\alpha}_1))$ and $\rightarrow I \xi^{\alpha}_3(\rightarrow I \xi^{\alpha}_1(\xi^{\alpha}_1))$ can be re-indexed to $\rightarrow I \xi^{\alpha}_4(\rightarrow I \xi^{\alpha}_1(\xi^{\alpha}_1))$.}

Finally, given two atomic bases $\mathfrak{B}_1, \mathfrak{B}_2$ on the same background language and with respective atomic systems $\texttt{S}_1, \texttt{S}_2$, and given a ground $g_1$ on $\mathfrak{B}_1$ and a ground $g_2$ on $\mathfrak{B}_2$, $g_1 \approx_{\mathfrak{B}_1/\mathfrak{B}_2} g_2$ can be defined in the same way as identity on one and the same atomic base, except that, in the case of atomic derivations, we require them to belong to $\texttt{S}_1$ and $\texttt{S}_2$ respectively, whereas in the second clause we replace $\approx_\mathfrak{B}$ with $\approx_{\mathfrak{B}_1/\mathfrak{B}_2}$ - and the first occurrence of $\mathfrak{B}$ with $\mathfrak{B}_1$, the second with $\mathfrak{B}_2$, $\approx_{\mathfrak{B}, \mathfrak{B}}$ clearly being $\approx_\mathfrak{B}$.

\section{Languages of grounding}

In this Section, we introduce languages of grounding and their expansions.

\subsection{Languages and expansions}

A language of grounding is relative to an atomic base. It contains operational symbols made of labels attached to operational types. Since we want languages of grounding to be ‘‘recursive", the set of the operational symbols is either finite, with a finite number of labels associated to an operational type $\tau$ in an equally finite set of operational types, or non-finite but partitioned after \emph{schemes of operational types} $\sigma$, where each scheme is associated to a finite number of labels. A scheme of operational types can be understood as an expression with meta-variables, indicating a structure that all its \emph{instances} are to have. We remark that we authorize only operational symbols with operational types such that operations on grounds of those types exist on the reference base.

\begin{defn}
Let $\mathfrak{B}$ be an atomic base with atomic system $\texttt{S}$, over a background language $L$, and to which a domain of individuals $D_\mathfrak{B}$ is associated. A \emph{language of grounding} $\Lambda$ over $\mathfrak{B}$ is specified by
\begin{itemize}
    \item an \emph{alphabet} $\texttt{Al}_\Lambda$ containing
    \begin{itemize}
        \item an individual constant $\delta_i$ that is a name of the $i$-th atomic derivation in $\texttt{S}$ with no undischarged assumptions and unbound individual variables;
        \item typed-variables $\xi^\alpha_i$ ($\alpha \in \texttt{FORM}_L$, $i \in \mathbb{N}$);
        \item operational symbols
        \begin{center}
            $F_1[\tau_i], ..., F_n[\tau_i]$ ($n \geq 1$)
        \end{center}
        for $\tau_i \in \{\tau_1, ..., \tau_m\}$ ($i \leq m$, $m \geq 0$) set of operational types of first-level over $L$, or
        \begin{center}
            $F_1^i[\tau_{\sigma_i}], ..., F_n^i[\tau_{\sigma_i}]$ ($n \geq 1$)
        \end{center}
        for every $\tau_{\sigma_i}$ instance of $\sigma_i \in \{\sigma_1, ..., \sigma_n\}$ ($i \leq n$, $n \geq 0$) set of schemes of operational types of first-level over $L$. Each operational symbol may bind a sequence of individual variables $\underline{x}$ or a sequence of typed-variables $\underline{\xi}$ - where $\underline{x}$ and $\underline{\xi}$, as well as bindings, have to be understood with the same notational conventions and with the same wording as in Section 4.2.
        
        We put the following restriction. If an operational type occurring in an operational symbol $\phi$ is $\alpha_1(\underline{x}_1), ..., \alpha_n(\underline{x}_n) \rhd \beta(\underline{x}_{n + 1})$, then there is a $\mathfrak{B}$-operation on grounds $f$ such that
        \begin{itemize}
            \item $f$ is defined on the same number of individuals as the length obtained by removing $\underline{x}$ from $\underline{x}_1, ..., \underline{x}_n$;
            \item the operational type of $f$ has co-domain $\beta(\underline{x}_{n + 1})$ and, for every $i \leq n$, the entry of index $i$ has co-domain $\alpha_i(\underline{x}_i)$ and domain the set of the types of the typed-variables that $\phi$ binds on $i$, if any, otherwise the entry of index $i$ is $\alpha_i(\underline{x}_i)$;
            \item $\phi$ binds $x$ and $\xi^\alpha$ on indexes $i$ and $j$ respectively iff $f$ binds $x$ and $\xi^\alpha$ on indexes $i$ and $j$ respectively ($i, j \leq n$);
        \end{itemize}
        \item a set of typed \emph{terms} $\texttt{TERM}_\Lambda$ specified by induction. If the operational type occurring in an operational symbol $\phi$ is $\alpha_1(\underline{x}_1), ..., \alpha_n(\underline{x}_n) \rhd \beta(\underline{x}_{n + 1})$, and if $\phi$ binds a sequence of individual variables $\underline{x}$ and a sequence of typed-variables $\underline{\xi}$, the inductive clause says - we indicate that a term has a certain type by colon -
            \begin{center} 
            $U_i : \alpha_i(\underline{x}_i) \in \texttt{TERM}_\Lambda \ (i \leq n) \Rightarrow \phi \ \underline{x} \ \underline{\xi} \ (U_1, ..., U_n) : \beta(\underline{x}_{n + 1}) \in \texttt{TERM}_\Lambda$.
            \end{center}
        The clause puts restrictions 1 and 2 as in Section 4.4 on the formation of the term.
    \end{itemize}
\end{itemize}
\end{defn}
\noindent Given a language of grounding $\Lambda$, let us indicate with $\mathbb{F}_\Lambda$ the set of its operational symbols.

\begin{defn}
Let $\Lambda_1$ be a language of grounding over an atomic base $\mathfrak{B}_1$. An \emph{expansion} of $\Lambda_1$ is language of grounding $\Lambda_2$ on an atomic base $\mathfrak{B}_2$ with $\mathfrak{B}_2$ expansion of $\mathfrak{B}_1$ and such that $\mathbb{F}_{\Lambda_1} \subseteq \mathbb{F}_{\Lambda_2}$.
\end{defn}

\noindent In line with what we have said about the non-monotonic character of the notion of operations on grounds, we remark that the expansion-relation between atomic bases may not coincide with the expansion-relation between languages of grounding. In other words, given a language of grounding $\Lambda$ on an atomic base $\mathfrak{B}_1$, there may be no expansion of $\Lambda$ on some expansion $\mathfrak{B}_2$ of $\mathfrak{B}_1$ - e.g., suppose that the alphabet of $\Lambda$ contains some operational symbol $F[\tau]$, which means that there is a $\mathfrak{B}_1$-operation on grounds of operational type $\tau$, but $\tau$ is no longer inhabited on $\mathfrak{B}_2$.

\subsection{Some examples}

\textbf{Example 1} Let $\mathfrak{B}$ be any atomic base with atomic system $\texttt{S}$ on a background language $L$. We call \emph{core language} on $\mathfrak{B}$ - indicated with $\texttt{G}$ - a language of grounding that, apart from names $\delta$ of derivations in $\texttt{S}$ with no undischarged assumptions and unbound individual variables, and apart from typed-variables $\xi^\alpha$, has the following operational symbols:
\begin{itemize}
    \item $\wedge I [\alpha, \beta \rhd \alpha \wedge \beta]$ ($\alpha, \beta \in \texttt{FORM}_L$);
    \item $\vee I [\alpha_i \rhd \alpha_1 \vee \alpha_2]$ ($i = 1, 2$, $\alpha_1, \alpha_2 \in \texttt{FORM}_L$);
    \item $\rightarrow I [\alpha \rhd \alpha \rightarrow \beta]$ binding $\xi^\alpha$ ($\alpha, \beta \in \texttt{FORM}_L$);
    \item $\exists I [\alpha(t) \rhd \exists x_i \alpha(x_i)]$ ($t \in \texttt{TERM}_L$, $\alpha(t/x_i) \in \texttt{FORM}_L$, $i \in \mathbb{N}$)
    \item $\forall I [\alpha(x_i) \rhd \forall x_j \alpha(x_j/x_i)]$ binding $x_i$ ($i, j \in \mathbb{N}$, $\alpha(x_i) \in \texttt{FORM}_L$);
    \item $\bot_\alpha [\bot \rhd \alpha]$ ($\alpha \in \texttt{FORM}_L$).
\end{itemize}
We qualify these operational symbols as \emph{primitive}. As an example of term formation, $\texttt{TERM}_{\texttt{G}}$ is specified in a standard inductive way like in the following cases - we omit type, subscripts and superscripts whenever possible:
\begin{itemize}
    \item $\delta : \alpha \in \texttt{TERM}_{\texttt{G}}$ - where $\alpha$ is the conclusion of the derivation of which $\delta$ is a name;
    \item $\xi^\alpha : \alpha \in \texttt{TERM}_{\texttt{G}}$;
    \item $T : \alpha_i \in \texttt{TERM}_{\texttt{G}} \Rightarrow \vee I [\alpha_i \rhd \alpha_1 \vee \alpha_2](T) : \alpha_1 \vee \alpha_2 \in \texttt{TERM}_{\texttt{G}}$;
    \item $T : \beta \in \texttt{TERM}_{\texttt{G}} \Rightarrow \ \rightarrow I \xi^\alpha (T) : \alpha \rightarrow \beta \in \texttt{TERM}_{\texttt{G}}$;
    \item $U : \alpha(x) \in \texttt{TERM}_{\texttt{G}} \Rightarrow \forall I [\alpha(x) \rhd \forall y \alpha(y/x)] x (U) : \forall y \alpha(y/x) \in \texttt{TERM}_{\texttt{G}}$ - with the restriction that $x \notin FV(\beta)$ for (see below for the notation) $\xi^\beta \in FV^T(U)$;
    \item $T : \bot \in \texttt{TERM}_{\texttt{G}} \Rightarrow \bot_\alpha(T) : \alpha \in \texttt{TERM}_{\texttt{G}}$.
\end{itemize}
Some notions must be taken as defined in a standard way on every language of grounding:

\begin{itemize}
    \item the set $S(T)$ of the \emph{sub-terms} of $T$;
    \item the sets $FV^I(T)$ and $BV^I(T)$ of, respectively, the \emph{free} and \emph{bound individual variables} of $T$;
    \item the sets $FV^T(U)$ and $BV^T(U)$ of, respectively, the \emph{free} and \emph{bound typed-variables} of $U$.
\end{itemize}
A term $U$ is said to be \emph{closed} iff $FV^I(U) = FV^T(U) = \emptyset$, \emph{open} otherwise. The following fact holds - and must hold for any language of grounding: given any operational symbol $\phi$ of $\texttt{G}$, for every $\phi \ \underline{x} \ \underline{\xi} \ (\dots U_i \dots) : \alpha \in \texttt{TERM}_{\texttt{G}}$,
\begin{enumerate}
    \item[a)] $FV^I(\phi \ \underline{x} \ \underline{\xi} \ (\dots U_i \dots))$ is equal to the union of $FV(\alpha)$ and of $FV^I(U_i)$ minus the individual variables in $\underline{x}$ with index $i$;
    \item[b)] $FV^T(\phi \ \underline{x} \ \underline{\xi} \ (\dots U_i \dots))$ is equal to $FV^T(U_i)$ minus the typed-variables in $\underline{\xi}$ with index $i$.
\end{enumerate}

\noindent \textbf{Example 2} A \emph{Gentzen-language} - indicated with $\texttt{Gen}$ - is an expansion of a core-language obtained by adding the \emph{non-primitive} operational symbols:
\begin{itemize}
    \item $\wedge_{E, i} [\alpha_1 \wedge \alpha_2 \rhd \alpha_i]$ ($i = 1, 2$, $\alpha_1, \alpha_2 \in \texttt{FORM}_L$);
    \item $\vee E [\alpha \vee \beta, \gamma, \gamma \rhd \gamma]$ - binding $\xi^\alpha$ on the second entry and $\xi^\beta$ on the third ($\alpha, \beta, \gamma \in \texttt{FORM}_L$);
    \item $\rightarrow E [\alpha \rightarrow \beta, \alpha \rhd \beta]$ ($\alpha, \beta \in \texttt{FORM}_L$);
    \item $\exists E [\exists x \alpha(x), \beta \rhd \beta]$ - binding $x$ and $\xi^{\alpha(x)}$ on the second entry ($\alpha(x), \beta \in \texttt{FORM}_L$);
    \item $\forall E [\forall x \alpha(x) \rhd \alpha(k)]$ ($\alpha(x) \in \texttt{FORM}_L$, $k \in \texttt{TERM}_L$ name of individual in $D_\mathfrak{B}$).
\end{itemize}
It is easily seen that $\texttt{Gen}$ is a kind of Curry-Howard linear translation of Gentzen's natural deduction system for first-order intuitionistic logic.
\medskip

\noindent \textbf{Example 3} Let $\mathfrak{B}$ be an atomic base with atomic system $\texttt{Ar}$ of Section 3.2.1 on a first-order language $L$ for arithmetic, to which the domain of individuals $\mathbb{N}$ is associated. A language of grounding for Heyting's arithmetic - indicated with $\texttt{HA}$ - is an expansion of $\texttt{Gen}$ over $\mathfrak{B}$ obtained by adding the following \emph{non-primitive} operational symbols:
\begin{itemize}
    \item $\texttt{Ind} [\alpha(0), \alpha(s(x)) \rhd \alpha(n)]$ - binding $x$ and $\xi^{\alpha(x)}$ on the second entry ($\alpha(0), \alpha(s(x)) \in \texttt{FORM}_L$, $n$ name of individual in $\mathbb{N}$)
\end{itemize}
Again, it is easily seen that $\texttt{HA}$ is a kind of Curry-Howard linear translation of a Gentzen's natural deduction system for Heyting's first-order arithmetic. We have of course G\"{o}del's incompleteness: there is closed $G \in \texttt{FORM}_L$ such that there is no closed $T : G \in \texttt{TERM}_{\texttt{HA}}$, and no closed $T : \neg G \in \texttt{TERM}_{\texttt{HA}}$.

\section{Denotation}

So far we have introduced formal languages which, as the examples at the end of the previous section show, can be understood as Curry-Howard translations of formal systems in a Gentzen's natural deduction style. But if our languages are to be languages \emph{of grounding}, we must explain how their terms ‘‘speak of" grounds and operations on grounds, i.e. how such terms denote the objects - described in Section 4 - based on which semantic notions are given. In order to do this, in this Section we introduce denotation functions. We first of all fix denotation for the elements of the alphabet, and then we extend this alphabetic denotation to terms in an inductive way.

The basic idea is that a closed term may denote a ground for a categorical assertion in a direct or indirect way, depending on whether it has a primitive operational symbol as outermost sign, or can be ‘‘computed" to such form. Computation rules are given by mapping the operational symbols of the term onto operations on grounds of the same type - as said in Section 4.3, operations on grounds are meant to come with defining equations that show how the operation is to be computed on privileged arguments. Observe that, if the type of the operational symbol is inhabited by \emph{two or more} distinct - not equivalent or not identical - operations on grounds, two or more distinct denotations for the same symbol are available. Open terms thus denote the operations on grounds that yield the grounds denoted by the closed instances of these terms, when applied to the grounds denoted by the closed terms through which the closed instances are obtained. Open terms are not computed, although a computational reducibility relation may hold between them - in this case they will denote equivalent, or even identical operations, depending on the denotation [a similar approach for Prawitz's valid arguments has been put forward in Tranchini 2014b, 2019 where, roughly, a non-normal derivation is taken to denote its so-called full-evaluated form].

\subsection{Denotation for alphabet and terms}

Below, as done with the notion of scheme of operational types, we employ the notion of \emph{scheme of} (\emph{systems of}) \emph{equations}, which must be understood as a meta-linguistic description of a class of equations sharing a certain structure.

\begin{defn}
Let $\Lambda$ be a language of grounding over an atomic base $\mathfrak{B}$ with atomic system $\texttt{S}$. A \emph{denotation for the alphabet} of $\Lambda$ is a function $den^*$ specified as follows:
\begin{itemize}
    \item $den^*(\delta_i) = i$-th derivation of $\texttt{S}$ of which $\delta_i$ is a name;
    \item $den^*(\xi^\alpha) = \texttt{Id}(FV(\alpha), \xi^\alpha)$, where $\texttt{Id}(FV(\alpha), \xi^\alpha)$ is the identity $\mathfrak{B}$-operation on grounds of type $\alpha \rhd \alpha$\footnote{There is thus a distinction between the typed-variables we used to speak about operations on grounds, and those in languages of grounding. We should have employed a distinct notation for the former, but we have preferred not to burden the notation.};
    \item given an operational symbol $\phi$,
    \begin{itemize}
        \item if $\phi$ is primitive, $den^*(\phi) = \phi$ - i.e. a primitive operational symbol corresponds to the primitive operation with the same label in the ground-clauses;
        \item if $\phi$ is non-primitive, $den^*(\phi)$ is one of the $\mathfrak{B}$-operations on grounds assumed as existing in connection with $\phi$ in definition 10.
        
        We have the following restriction: for every two operational symbols $F[\tau^1_\sigma]$ and $F[\tau^2_\sigma]$ where $\tau^1_\sigma$ and $\tau^2_\sigma$ are instances of the same scheme of operational types, $den^*(F[\tau^1_\sigma])$ and $den^*(F[\tau^2_\sigma])$ are $\mathfrak{B}$-operations on grounds whose defining (systems of) equations are instances of the same scheme of (systems of) equations. 
    \end{itemize}
\end{itemize}
\end{defn}

\begin{defn}
Let $\Lambda$ be a language of grounding, and let $den^*$ be a denotation for the elements of its alphabet. A \emph{denotation for the terms} of $\Lambda$ \emph{associated to} $den^*$ is a function $den$ specified as follows:
\begin{itemize}
    \item $den(\delta) = den^*(\delta)$;
    \item $den(\xi^\alpha) = den^*(\xi^\alpha)$;
    \item $den(\phi \ \underline{x} \ \underline{\xi} \ (U_1, ..., U_n)) = den^*(\phi)(\underline{y}, den(U_1), ..., den(U_n))$.\footnote{$\underline{y}$ is the sequence of the individual variables of the operational type of $den^*(\phi)$ if $den^*(\phi)$ is an operation on grounds of first-level, and those of the co-domains of the operational type of $den^*(\phi)$ if $den^*(\phi)$ is an operation on grounds of second-level.}
\end{itemize}
\end{defn}
\noindent As an example, take the languages of grounding $\texttt{Gen}$ and $\texttt{HA}$ in Section 5.2, and consider the defining equations of Section 4.3. We set $den^*(k E) = f_k$ for each logical constant $k$, and $den^*(\texttt{Ind}[\alpha(0), \alpha(s(x)) \rhd \alpha(n)]) = f^n_\texttt{Ind}$.

The non-monotonic character of operations on grounds also affects the possibility of expanding denotation functions - similarly to what happens with expansions of atomic bases and expansions of languages of grounding. Given a language of grounding $\Lambda$ on an atomic base $\mathfrak{B}_1$, and a denotation function $den^*$ for the elements of the alphabet of $\Lambda$, let $\mathfrak{B}_2$ be an expansion of $\mathfrak{B}_1$. Two cases may occur: (1) either, as already seen at the end of Section 5.1, there is no expansion of $\Lambda$ over $\mathfrak{B}_2$, so that $den^*$ cannot be extended to $\mathfrak{B}_2$, or (2) there is an expansion of $\Lambda$ over $\mathfrak{B}_2$, but $den^*$ associates to some operational symbol of $\Lambda$ some $\mathfrak{B}_1$-operation on grounds fixed by an equation which does not also define a $\mathfrak{B}_2$-operation on grounds of the same operational type. To illustrate the second case, suppose that $\mathfrak{B}_1$ is a logical base on a language for first-order arithmetic, that the alphabet of $\Lambda$ contains just one operational symbol of operational type $0 = 0 \rhd 1 = 1$, and finally that $den^*$ associates to this symbol the empty function of type $0 = 0 \rhd 1 = 1$ - i.e. the function with null set of computation instructions. Let now $\mathfrak{B}_2$ be the expansion of $\mathfrak{B}_1$ whose atomic system is $\texttt{Ar}$ of Section 3.2.1. Clearly, the operational type $0 = 0 \rhd 1 = 1$ is still inhabited on $\mathfrak{B}_2$ - so $\Lambda$ can be expanded on the new base - but it no longer contains the empty function, so $den^*$ cannot be defined on the expansion of $\Lambda$ on $\mathfrak{B}_2$.

\subsection{Denotation theorem}

The following denotation theorem is a kind of ‘‘correctness" result for languages of grounding, as it tells us that any term of any such language denotes a ground on the reference base. The result obviously depends on the restriction we put in definition 10 on operational symbols, i.e. the fact that we allowed only operational symbols whose operational type is inhabited by an adequate operation on grounds.

\begin{thm}
Let $\Lambda$ be a language of grounding on an atomic base $\mathfrak{B}$ to which a domain of individuals $D_\mathfrak{B}$ is associated, let $den^*$ be a denotation function for the elements of $\Lambda$, and let $den$ be a denotation function for the terms of $\Lambda$ associated to $den^*$. Let $U : \beta \in \texttt{TERM}_\Lambda$ be such that
\begin{center}
    $FV^I(U) = \{x_1, ..., x_n\}$ and $FV^T(U) = \{\xi^{\alpha_1}, ..., \xi^{\alpha_m}\}$,
\end{center}
and let $\{\alpha_s, ..., \alpha_t\}$ be the set of the formulas that occur as types of the elements of $FV^T(U)$. $den(U)$ is a ground on $\mathfrak{B}$ for $\alpha_s, ..., \alpha_t \vdash \beta$ defined on $n$ individuals on $D_\mathfrak{B}$.
\end{thm}

\begin{proof}
By simple induction on the complexity of $U$, although there are many cases, each of which has in turn many sub-cases. One has to use facts a) and b) in Section 5.2.
\end{proof}
\noindent The converse of theorem 14 would accordingly be a kind of ‘‘completeness" result, saying that every ground on an atomic base is denoted by some term in some language of grounding on that base, i.e.

\begin{itemize}
    \item[(K)] for every atomic base $\mathfrak{B}$ and $g$ ground on $\mathfrak{B}$ for $\alpha_1, ..., \alpha_n \vdash \beta$ ($n \geq 0$), there is a language of grounding $\Lambda$ on $\mathfrak{B}$ such that, for some denotation function $den^*$ for the elements of the alphabet of $\Lambda$, and called $den$ the denotation function for the terms of $\Lambda$ associated to $den^*$, there is $U : \beta \in \texttt{TERM}_\Lambda$ such that $den(U) \approx_\mathfrak{B} g$.
\end{itemize}

\section{Internal properties}

We qualify as \emph{internal} the first kind of properties of a language of grounding we discuss. These properties are called internal because they concern a language of grounding ‘‘in itself", namely, they are properties of the set of the terms of the language.

\subsection{Closure under canonical form}

\begin{defn}
Let $\Lambda$ be a language of grounding and $T \in \texttt{TERM}_\Lambda$. We say that $T$ is \emph{canonical} iff it is an individual constant, a typed-variable, or its outermost operational symbol is one of the operational symbols of a core-language of grounding, namely $\wedge I$, $\vee I$, $\rightarrow I$, $\exists I$ or $\forall I$. $T$ is \emph{non-canonical} otherwise.
\end{defn}
\noindent Now, the denotation theorem tells us that, whatever the denotation on $\Lambda$ is, if $T$ is closed, it denotes a ground for a categorical judgment or assertion, i.e. an object whose outermost operation is primitive. We may therefore expect that it is possible to find a canonical term in $\Lambda$ having a denotation equivalent to that of $T$, or even the same denotation as $T$ - according to whether we consider equivalence or identity as introduced in Section 4.5.

In the following definition, as well as on some consequences we draw from it, we focus on equivalence, but identity is discussed immediately afterwards. We just remark that, from now on, we shall leave implicit that we are working on a language of grounding $\Lambda$ (possibly with indexes) over an atomic base $\mathfrak{B}$ (with atomic system $\texttt{S}$ and over a background language $L$), with $den^*$ (possibly with indexes) denotation function for the alphabet of $\Lambda$ and $den$ (possibly with indexes) denotation function for the terms of $\Lambda$ associated to $den^*$. This is because we prefer not to burden the notation excessively.
\begin{defn}
$\Lambda$ is \emph{closed under canonical form relative to} $den$ iff, for every closed $T \in \texttt{TERM}_\Lambda$, there is a closed canonical $U \in \texttt{TERM}_\Lambda$ such that $den(T) \equiv_\mathfrak{B} den(U)$.
\end{defn}
We have of course examples of languages of grounding closed under canonical form relative to a denotation defined on their terms - e.g. $\texttt{Gen}$, via the denotation for the alphabet suggested in Section 6.1, and thanks to Prawitz's normalisation theorems [Prawitz 2006]. But in general, this does not hold for every language of grounding - think of a language of grounding that is a kind of Curry-Howard translation of a system the closed derivations of which cannot always be reduced to canonical form in the same system.

Nonetheless, we may ask if any language of grounding can be at least expanded to one that enjoys such property relative to some denotation. Concerning this point, we prove a weaker result, namely that, if (K) in Section 6.2 holds, any language of grounding can be expanded to one that is closed under canonical form \emph{with respect to a given closed term}. We leave aside whether the expansion can be found also for the whole set of closed terms.
\begin{defn}
Let $T \in \texttt{TERM}_\Lambda$ be closed. We say that $\Lambda$ is \emph{closed under canonical form relative to} $den$ \emph{with respect to} $T$ iff, there is a closed canonical $U \in \texttt{TERM}_\Lambda$ such that $den(T) \equiv_\mathfrak{B} den(U)$.
\end{defn}
\begin{thm}
Suppose that (K) in Section 6.2 holds. Then, given $\Lambda_1$ over $\mathfrak{B}$ with atomic system $\texttt{S}$, $den^*_1$ over $\Lambda_1$ and $T \in \texttt{TERM}_{\Lambda_1}$ closed, there is an expansion $\Lambda_2$ of $\Lambda_1$ on $\mathfrak{B}$ such that, for some $den^*_2$ over $\Lambda_2$, $den_1(T) \approx_\mathfrak{B} den_2(T)$ and $\Lambda_2$ is closed under canonical form relative to $den_2$ with respect to $T$.
\end{thm}
\begin{proof}
By cases on the type $\alpha$ of $T$ - for $\alpha$ of course closed. If $\alpha$ is atomic, then $den_1(T)$ is an atomic derivation $\Delta$ in $\texttt{S}$ named by some $\delta$ in $\Lambda_1$. We can hence put $\Lambda_2 = \Lambda_1$ and $den_2 = den_1$, since
\begin{center}
    $den_1(T) = \Delta = den_1(\delta)$
\end{center}
and $\Delta \equiv_\mathfrak{B} \Delta$. For $\alpha$ logically complex, we consider the case $\alpha = \beta_1 \vee \beta_2$ - for $\beta_i$ of course closed ($i = 1, 2$). $den_1(T)$ is in this case $\vee I [\beta_i \rhd \beta_1 \vee \beta_2](g)$ for $g$ ground on $\mathfrak{B}$ for $\vdash \beta_i$. Since (K) holds, there is $\Lambda^*$ on $\mathfrak{B}$ such that, for some closed $U : \beta_i \in \texttt{TERM}_{\Lambda^*}$, for some $den^*$, $den(U) \approx_\mathfrak{B} g$. The only difficult case is the following. $\texttt{Al}_{\Lambda_1} \cap \ \texttt{Al}_{\Lambda^*} \neq \emptyset$ and, by indicating with $\mathbb{F}_{\Lambda_1, T}$ and $\mathbb{F}_{\Lambda^*, U}$ the sets of the non-primitive operational symbols of $\Lambda_1$ and $\Lambda^*$ respectively occurring in $T$ and $U$ respectively, $\mathbb{F}_{\Lambda_1, T} \cap \mathbb{F}_{\Lambda^*, U} \neq \emptyset$ and, for some $x \in \mathbb{F}_{\Lambda_1, T} \cap \mathbb{F}_{\Lambda^*, U}$, it does not hold that $den^*_1(x) \approx_\mathfrak{B} den^*(x)$. For every operational symbol $F[\tau] \in \mathbb{F}_{\Lambda_1, T} \cap \mathbb{F}_{\Lambda^*, U}$ for which it does not hold that $den^*_1(F[\tau]) \approx_\mathfrak{B} den^*(F[\tau])$, let us take an operational symbol $F_\dagger[\tau]$ such that $F_\dagger[\tau] \notin \mathbb{F}_{\Lambda_1} \cup \mathbb{F}_{\Lambda^*}$. Let us call $\mathbb{F}_\dagger$ the set of the operational symbols we obtain in this way. We let $\Lambda_2$ be such that $\texttt{Al}_{\Lambda_2} = (\texttt{Al}_{\Lambda_1} \cup \ \texttt{Al}_{\Lambda^*}) \cup \mathbb{F}_\dagger$, and $den^*_2$ be such that
    \begin{center}
        $den^*_2(x) = \begin{cases} den^*_1(x) & \text{if} \ x \in \texttt{Al}_{\Lambda_1} \\ den^*(x) & \text{if} \ x \in (\texttt{Al}_{\Lambda^*} - (\texttt{Al}_{\Lambda_1} \cap \ \texttt{Al}_{\Lambda^*})) \\ den^*(F[\tau]) & \text{if} \ x = F_\dagger[\tau] \end{cases}$
    \end{center}
Let $U_\dagger$ be the term of $\texttt{TERM}_{\Lambda_2}$ obtained from $U$ by replacing with $F_\dagger[\tau] \in \mathbb{F}_\dagger$ any $F[\tau] \in \mathbb{F}_{\Lambda_1, T} \cap \mathbb{F}_{\Lambda^*, U}$ for which it does not hold that $den^*_1(F[\tau]) \approx_\mathfrak{B} den^*(F[\tau])$. We have that $den_2(T) \approx_\mathfrak{B} den_1(T)$ and $den_2(U_\dagger) \approx_\mathfrak{B} den(U)$, so that

\begin{center}
     $den_1(T) = \vee I [\beta_i \rhd \beta_1 \vee \beta_2](g) = den^*(\vee I [\beta_i \rhd \beta_1 \vee \beta_2])(den^*(U)) = den(\vee I [\beta_i \rhd \beta_1 \vee \beta_2] (U))$
\end{center}
and clearly $\vee I [\beta_i \rhd \beta_1 \vee \beta_2](g) \equiv _\mathfrak{B} \vee I [\beta_i \rhd \beta_1 \vee \beta_2](g)$.
\end{proof}
Observe that theorem 18 holds also in a weaker form, i.e. by replacing $\approx_\mathfrak{B}$ with $\equiv_\mathfrak{B}$ in (K). Likewise, what we have said so far for closure under canonical form can be also given in terms of identity on $\mathfrak{B}$. It is sufficient to replace $\equiv_\mathfrak{B}$ with $\approx_\mathfrak{B}$ in definition 16. This notion is somewhat more stringent than the one identified by definition 16, so we may speak in this case of $\Lambda$ being \emph{strictly closed under canonical form}. Then, an example analogous to the one provided for simple closure under canonical form, where $\equiv_\mathfrak{B}$ is replaced with $\approx_\mathfrak{B}$, shows that not every language of grounding is strictly closed under canonical form. By replacing again $\equiv_\mathfrak{B}$ with $\approx_\mathfrak{B}$ in definition 17 we obtain an analogous notion of \emph{strictly closed under canonical form with respect to a closed term}. Theorem 18 still holds - indeed, our proof \emph{actually} is for the stricter notion.

We finally remark that the strategy adopted for proving theorem 18 is based on what we may call a ‘‘duplication" of an operational symbol $F[\tau]$ occurring both in $T$ and in $U$, i.e. we replace the label $F$ with a fresh label $F_\dagger$. Now, if operational symbols are understood as inferences (instantiating given inferences rules), this strategy would make little sense, as one usually identifies an inference simply by its premises and its conclusion - so, the operational type $\tau$ being equal, we are not allowed to consider $F[\tau]$ and $F_\dagger[\tau]$ as distinct operational symbols. But we should not forget that we are considering operational symbols \emph{interpreted} by a denotation function, which would correspond to an inference (instantiating a given inference rule) \emph{equipped with} a justification procedure for this inference (rule) as in Prawitz's semantics of valid arguments.

\subsection{Universal denotation}

So far, we have spoken of grounds \emph{relative to an atomic base}. Likewise, languages of grounding have been defined \emph{over an atomic base}. One may naturally wonder whether we can speak of grounds over an atomic base that remain such \emph{whatever the atomic base is}, and what kind of terms in languages of grounding would denote these ‘‘universal" grounds.

In the case of grounds, we proceed as follows. Given an atomic base $\mathfrak{B}_1$ on a background language $L_1$, and a ground $g_1$ on $\mathfrak{B}_1$, $g_1$ is \emph{universal} iff, for every $L_2$ expansion of $L_1$, for every atomic base $\mathfrak{B}_2$ on $L_2$, there is a ground $g_2$ on $\mathfrak{B}_2$ such that $g_1 \approx_{\mathfrak{B}_1, \mathfrak{B}_2} g_2$ - as the notion of identity over two bases has been introduced in Section 4.5. Observe that $L_1$ and $L_2$ are here \emph{background languages}, i.e. languages for the typing of $g_1$ and $g_2$, not languages of grounding and expansions of these languages as given in definitions 10 and 11. Therefore, $\mathfrak{B}_2$ is not required to be an expansion of $\mathfrak{B}_1$.

Universality as defined here corresponds to what, in Prawitz's semantics of valid arguments, is sometimes called \emph{schematic validity}. A ground $g_1$ on $\mathfrak{B}_1$ is universal when, for every $\mathfrak{B}_2$ on expansions of the background language of $\mathfrak{B}_1$, there is a ground $g_2$ on $\mathfrak{B}_2$ that is ‘‘computationally" equal to $g_1$. It is therefore not sufficient that $g_2$ be a ground on $\mathfrak{B}_2$ for the same assertion as the one which $g_1$ is a ground for on $\mathfrak{B}_1$; $g_2$ must more strongly be the ‘‘same" ground as $g_1$, i.e. $g_1$ and $g_2$ must involve operations defined by the same schemes of equations, and composed according to the same order. Thus, for it to be universal, $g_1$ must be the ‘‘same" ground everywhere. Of course, we cannot expect that the operations involved in $g_1$ have everywhere the same domain and co-domain - objects within domain and co-domain will vary depending on the base; this is why we defined universality in terms of identity, i.e. of equality of computation instructions throughout the bases. The idea is that the operations involved in $g_1$ can be defined in the same way over all bases, without this altering the fact that $g_1$ is a ground on the base. E.g., consider $f_\wedge$ of operational type $\alpha_1 \wedge \alpha_2 \rhd \alpha_i$ ($i = 1, 2$) defined as in Section 4.3 relative to an arbitrary base $\mathfrak{B}$, i.e. for every $g_1$ ground on $\mathfrak{B}$ for $\vdash \alpha_1$, for every $g_2$ ground on $\mathfrak{B}$ for $\vdash \alpha_2$,

\begin{center}
    $f_\wedge(\wedge I (g_1, g_2)) = g_i$.
\end{center}
One easily proves that this is a ground on $\mathfrak{B}$ for $\alpha_1 \wedge \alpha_2 \vdash \alpha_i$, but since $\mathfrak{B}$ is arbitrary, the \emph{same equation} works everywhere, so that $f_\wedge$ is a universal ground for $\alpha_1 \wedge \alpha_2 \vdash \alpha_i$. The same can be proved for all the other operations defined in Section 4.3, except induction. Now, compare this to the following. Suppose that, for every atomic base $\mathfrak{B}_1$, there is a set of computation instructions $\varepsilon_{\mathfrak{B}_1}$ that defines a $\mathfrak{B}_1$-operation on grounds of operational type, say, $\alpha \rhd \beta$. Suppose further that, for any such $\varepsilon_{\mathfrak{B}_1}$, there is a base $\mathfrak{B}_2$ such that the same set of computation instructions \emph{does not} define a $\mathfrak{B}_2$-operation on grounds of operational type $\alpha \rhd \beta$. Thus, on every atomic base $\mathfrak{B}$ we have a ground for $\alpha \vdash \beta$ thanks to the $\mathfrak{B}$-operation on grounds defined by $\varepsilon_\mathfrak{B}$, but since every set of computation instructions fails on some bases, this ground will not be the ‘‘same" on all bases, so none of these grounds will be universal. This also concerns a difference between two ways of defining logical validity, which we shall refer to in Sections 9.1 and 9.4. We now move to universality of terms, for which we give the following definition.

\begin{defn}
Given $\Lambda$ and $den^*$, $T \in \texttt{TERM}_\Lambda$ is \emph{universal with respect to} $den^*$ iff, for every $x \in \texttt{Al}_\Lambda$ occurring in $T$, $den^*(x)$ is a universal ground.
\end{defn}

\begin{thm}
Given $\Lambda$ over $\mathfrak{B}$ and $den^*$, $T \in \texttt{TERM}_\Lambda$ is universal with respect to $den^*$ iff
\begin{itemize}
    \item[a)] $S(T)$ does not contain individual constants and
    \item[b)] for every non-primitive operational symbol $\phi$ of $\Lambda$ occurring in $T$, $den^*(\phi)$ is a universal ground.
\end{itemize}
\end{thm}
\begin{proof}
Straightforward.
\end{proof}
\begin{crl}
$T \in \texttt{TERM}_\Lambda$ is universal with respect to $den^*$ iff, for every $U \in S(T)$, $U$ is universal with respect to $den^*$.
\end{crl}

\noindent By assuming a plausible convention, universality can be extended to denotation functions for terms.
\begin{cnv}
Let $\mathfrak{B}$ be an atomic base and let $f(\underline{x}, \xi^{\tau_1}, ..., \xi^{\tau_n})$ be a $\mathfrak{B}$-operation on grounds that is, more in particular, a universal ground. Then, the operation is defined in such a way that, for every sequence of individuals $\underline{k}$ whose length is less-equal than that of $\underline{x}$, for every $g_i$ ground on $\mathfrak{B}$ for $\tau_i[\underline{k}/\underline{x}]$ ($i \leq n$, but not necessarily for all such $i$-s) that is more in particular a universal ground, $f(\underline{k}/\underline{x}, \dots, g_i, \dots)$ is a universal ground.
\end{cnv}
\begin{prp}
Given $\Lambda$ and $den^*$, and $T \in \texttt{TERM}_\Lambda$ universal with respect to $den^*$, $den(T)$ is a universal ground.
\end{prp}
\begin{proof}
By induction on the complexity of $T$, using corollary 21 and convention 22.
\end{proof}
\noindent The converse of proposition 23, however, does not hold. This is simply because a non-canonical term may ‘‘reduce" to a universal ground, and nonetheless involve some sub-term which is not universal - think of a non-canonical term of type $\alpha$ obtained by eliminating conjunction on a canonical term of type $\alpha \wedge \beta$, the immediate sub-term of type $\beta$ of which is not universal.

\section{External properties}

We finally discuss what we shall call \emph{external} properties of a language of grounding. At variance with the internal ones, they concern a relation that a language of grounding may or not have with other languages of grounding.

\begin{defn}
Given $\Lambda_1$ over $\mathfrak{B}_1$ with atomic system $\langle L_1, \Re_1 \rangle$, let $\Lambda_2$ be an expansion of $\Lambda_1$ over $\mathfrak{B}_2$ with atomic system $\langle L_2, \Re_2 \rangle$. We call $\Lambda_2$ a \emph{primitive} expansion of $\Lambda_1$ iff $\Re_1 \subset \Re_2$.
\end{defn}
\noindent Thus, $\Lambda_2$ primitively expands $\Lambda_1$ simply when the set of the individual constants of $\Lambda_1$ is a proper subset of the set of the individual constants of $\Lambda_2$.

There is of course another type of primitive expansion, namely, an expansion whose base has a background language with new \emph{constants}. The new constants, in fact, require the introduction of new primitive operations on grounds in terms of which the meaning of these constants can be determined, and correspondingly, new primitive operational symbols in the languages. Although this is surely true, we shall not investigate this possibility here - we just suggest two examples at the end of this section.

\begin{defn}
Given $\Lambda_1$ over an atomic base with background language $L$, let $\Lambda_2$ be an expansion of $\Lambda_1$ over an atomic base $\mathfrak{B}$, and let $den^*$ be relative to $\Lambda_2$. $\Lambda_2$ is \emph{conservative over} $\Lambda_1$ \emph{with respect to} $den^*$ iff, for every $U : \beta \in \texttt{TERM}_{\Lambda_2}$ with $\beta \in \texttt{FORM}_L$, and such that, for every $\xi^\alpha \in FV^T(U)$, $\alpha \in \texttt{FORM}_L$, there is $Z \in \texttt{TERM}_{\Lambda_1}$ such that $den(U) \equiv_\mathfrak{B} den(Z)$.
\end{defn}
\noindent Thus, intuitively, $\Lambda_2$ conservatively expands $\Lambda_1$ when, for every ground on $\mathfrak{B}$ for a judgment or assertion $\tau_1, ..., \tau_n \vdash \tau_{n + 1}$ on $L$ denoted via $den$ by some term of $\Lambda_2$, there is a term in $\Lambda_1$ that denotes via $den$ an equivalent ground. This means in particular that, for every operation whose operational type is on $L$, and which is denoted by some term of $\Lambda_2$, there is an operation denoted by some term of $\Lambda_1$ such that the two operations are ‘‘extensionally" equal - although they may computationally differ. So, all that can be expressed through $\Lambda_2$ on $\mathfrak{B}$, can be expressed through $\Lambda_1$ in an ‘‘extensionally" equal way.

The kind of conservativity we are dealing with here differs from the usual notion of conservativity of theories. We could speak of conservativity ‘‘of denotation", to distinguish it from conservativity of derivability. Between the two notions, however, a link obviously exists, which is established by the following theorem.
\begin{thm}
Given $\Lambda_1$ with background language $L$, let $\Lambda_2$ be an expansion of it over $\mathfrak{B}$, and let $den^*$ be relative to $\Lambda_2$. Let us finally suppose that, for every non-primitive operational symbol $\phi$ of $\Lambda_2$ whose operational type is on $L$, there is a composition of $\mathfrak{B}$-operations on grounds $f_1, ..., f_n$ resulting in a ground $h$ on $\mathfrak{B}$ such that:
\begin{itemize}
    \item for every $i \leq n$, there is an operational symbol $\phi_i$ of $\Lambda_1$ such that $den^*(\phi_i) \equiv_\mathfrak{B} f_i$, and
    \item $den^*(\phi) \equiv_\mathfrak{B} h$.
\end{itemize}
Then, $\Lambda_2$ is conservative over $\Lambda_1$ with respect to $den^*$.
\end{thm}
\begin{proof}
By induction on the complexity of $T$.
\end{proof}

So far, conservativity has been dealt with in terms of equivalence. We quickly define it as relative to identity.
\begin{defn}
Let $\Lambda_1$ and $\Lambda_2$ be as in definition 25 on $\mathfrak{B}_1$ and $\mathfrak{B}_2$ respectively. $\Lambda_2$ is \emph{strictly conservative over} $\Lambda_1$ \emph{with respect to} $den^*$ iff, for every $U : \beta \in \texttt{TERM}_{\Lambda_2}$ with $\beta \in \texttt{FORM}_L$ and such that, for every $\xi^{\alpha} \in FV^T(U)$, $\alpha \in \texttt{FORM}_L$, there is $Z \in \texttt{TERM}_{\Lambda_1}$ such that $den(U) \equiv_{\mathfrak{B}_2} den(Z)$ and, for some ground $g$ on $\mathfrak{B}_1$, $den(Z) \approx_{\mathfrak{B}_1, \mathfrak{B}_2} g$.
\end{defn}

\noindent We conclude with some examples showing that primitiveness and conservativity are not equivalent notions.
\medskip

\noindent \textbf{Non-primitive non-conservative} $\texttt{Gen}$ is a non-primitive expansion of a core-language $\texttt{G}$ when it has the same base as $\texttt{G}$. However, for no denotation function $\texttt{Gen}$ can be conservative over $\texttt{G}$, as in $\texttt{G}$ we have no term of type, e.g., $\alpha_1 \wedge \alpha_2 \rightarrow \alpha_i$ ($i = 1, 2$).
\medskip

\noindent \textbf{Non-primitive conservative} In Section 8.2 we have seen that a non-primitive expansion obtained by adding to $\texttt{Gen}$ an operational symbol $\phi$ that denotes, via $den^*$, $DS$ in Section 4.5 is conservative over $\texttt{Gen}$ with respect to $den^*$ when the latter assigns to the non-primitive operational symbols of $\texttt{Gen}$ the operations on grounds illustrated in Section 4.3.
\medskip

\noindent \textbf{Primitive non-conservative} Any expansion obtained by enlarging the atomic base - i.e. by adding individual constants. Alternatively, an expansion obtained by adding a reflection principle to $\texttt{HA}$ as in section 5.2
\medskip

\noindent \textbf{Primitive conservative} Consider $\Lambda_2$ obtained by adding to $\Lambda_1$ new logical constants. Suppose that $\Lambda_1$ and $\Lambda_2$ are a kind of Curry-Howard translations of respective Gentzen's natural deduction systems $\Sigma_1$ and $\Sigma_2$. Suppose finally that to the elimination rules of $\Sigma_2$ we can associate reduction procedures permitting to prove a sub-formula property. So, any derivation  in $\Sigma_2$ with assumptions and conclusion in the language of $\Sigma_1$ can be reduced to a derivation in $\Sigma_1$. And since we can assume that reduction procedures can be Curry-Howard translated into defining (systems of) equations for non-primitive operations to be associated to the non-primitive operational symbols of $\Lambda_2$, we can define on $\Lambda_2$ a $den^*$ in such a way that $\Lambda_2$ is conservative over $\Lambda_1$ with respect to $den^*$.

\section{Correctness, completeness and recognisability}

In this section we deal with correctness and completeness issues for intuitionistic first-order logic, based on a ground-theoretic definition of the notion of inferential validity.

\subsection{Inferential validity and correctness of intuitionistic logic}

As said in the beginning of this paper, the main aim of ToG is that of defining the notion of valid inference. This is done by requiring that an inference is valid when an operation on grounds exists from grounds for its premises to grounds for its conclusion, and with adequate bindings.
\begin{defn}
Let $\mathfrak{B}$ be an atomic base on a background language $L$. An inference on $L$
\begin{prooftree}
\AxiomC{$\alpha_1, ..., \alpha_n$}
\RightLabel{$R$}
\UnaryInfC{$\beta$}
\end{prooftree}
- possibly binding individual variables and discharging assumptions - is \emph{valid over} $\mathfrak{B}$ iff there is a $\mathfrak{B}$-operation on grounds $f$ of operational type $\tau_1, ..., \tau_n \rhd \beta$ where:
\begin{itemize}
    \item $\tau_i$ is $\alpha_i$ if $R$ does not discharge assumptions on $\alpha_i$, otherwise it is $\Gamma \rhd \alpha_i$, where $\Gamma$ is the set of formulas discharged by $R$ on $\alpha_i$ ($i \leq n$);
    \item  $f$ binds $x$ and $\xi^\gamma$ on index $i$ iff $R$ binds $x$ and discharges $\gamma$ on $\alpha_i$;
    \item $f$ is a proper ground for $\tau_1, ..., \tau_n \vdash \beta$.
\end{itemize}
We say that $R$ is \emph{logically valid} iff $f$ is a universal ground.
\end{defn}
\noindent A $\mathfrak{B}$-proof of $\alpha$ from $\Gamma$ is a chain of inferences valid on $\mathfrak{B}$. When such a proof exists, $\alpha$ is said to be a $\mathfrak{B}$-consequence of $\Gamma$ - written $\Gamma \models_\mathfrak{B} \alpha$; $\alpha$ is a logical consequence of $\Gamma$ - written $\Gamma \models \alpha$ - when there is a proof of $\alpha$ from $\Gamma$ that only involves logically valid inferences. It is immediate to observe that $\Gamma \models_\mathfrak{B} \alpha$ (resp. $\Gamma \models \alpha$) iff there is a composition of $\mathfrak{B}$-operation on grounds (resp. universal operations on grounds) resulting in a $\mathfrak{B}$-operation on grounds (resp. universal operation on grounds) of operational type $\Gamma \rhd \alpha$. With this established, it is very easy to prove correctness of first order intuitionistic logic - written $\texttt{IL}$.

\begin{thm}[Correctness of $\texttt{IL}$]
$\Gamma \vdash_{\texttt{IL}} \alpha \Rightarrow \Gamma \models \alpha$.
\end{thm}
\begin{proof}
The inferences of $\texttt{IL}$ are logically valid via the universal operations on grounds - except induction - defined in Section 5.2.
\end{proof}

\subsection{(In)completeness issues for intuitionistic logic}

\noindent The question whether $\texttt{IL}$ is also complete with respect to ToG amounts to so-called Prawitz's conjecture (Prawitz 1973). In our framework, the conjecture runs as follows.
\begin{cnj}[Prawitz's conjecture]
Given an operational type $\tau_1, ..., \tau_n \rhd \tau_{n + 1}$, if there is a universal operation on grounds $f$ of such type, then the inference
\begin{prooftree}
\AxiomC{$\sigma_1, ..., \sigma_n$}
\RightLabel{$R$}
\UnaryInfC{$\beta$}
\end{prooftree}
is derivable in $\texttt{IL}$, where:
\begin{itemize}
    \item $\sigma_i$ is $\tau_i$ if $\tau_i$ has empty domain, and
    \begin{prooftree}
    \AxiomC{$\Gamma_i$}
    \noLine
    \UnaryInfC{$\vdots$}
    \noLine
    \UnaryInfC{$\alpha_i$}
    \end{prooftree}
    if $\tau_i$ has domain $\Gamma_i$ and co-domain $\alpha_i$ ($i \leq n$);
    \item $f$ binds $x$ and $\xi^\gamma$ on index $i$ iff $R$ binds $x$ and discharges $\gamma$ on $\sigma_i$ ($i \leq n$);
    \item $\beta$ is the co-domain of $\tau_{n + 1}$.
\end{itemize}
\end{cnj}

\begin{prp}
Prawitz's conjecture implies that $\Gamma \models \alpha \Rightarrow \Gamma \vdash_{\texttt{IL}} \alpha$.
\end{prp}

Conjecture 30 requires that, if a scheme of operational types is ‘‘inhabited" by some universal operation on grounds with certain bindings, then the operational type is derivable in $\texttt{IL}$ with appropriate bindings of individual variables and discharging of assumptions. But the type could be ‘‘inhabited" by ‘‘extensionally" different universal operations on grounds. On the other hand, derivations in $\texttt{IL}$ correspond as seen to terms of $\texttt{Gen}$. Therefore, we may ask whether \emph{each} of the universal operations ‘‘inhabiting" the type is extensionally equal to the denotation of some term of $\texttt{Gen}$. This would amount to a stronger full-completeness conjecture.

\begin{cnj}[Full-completeness conjecture]
Given $\Lambda$ over $\mathfrak{B}_1$, let $den^*_1$ be relative to $\Lambda$ and such that, for some operational symbol $\phi$ of $\Lambda$, $den^*_1(\phi)$ is a universal operation on grounds $f$. Then, for some $den^*_2$ relative to $\texttt{Gen}$ over $\mathfrak{B}_1$, for some $\phi_1, ..., \phi_n$ operational symbols of $\texttt{Gen}$, there is a universal operation on grounds $h$ composed of $den^*_2(\phi_1), ..., den^*_2(\phi_n)$ such that, for every $\mathfrak{B}_2$, called $g_1$ the ground on $\mathfrak{B}_2$ such that $den^*_1(\phi) \approx_{\mathfrak{B}_1, \mathfrak{B}_2} g_1$, and called $g_2$ the ground on $\mathfrak{B}_2$ such that $h \approx_{\mathfrak{B}_1, \mathfrak{B}_2} g_2$, it holds that $g_1 \equiv_{\mathfrak{B}_2} g_2$.
\end{cnj}
\noindent In the case when $\phi$ does not bind individual or typed-variables, conjecture 32 boils down to requiring the existence of some $T \in \texttt{TERM}_{\texttt{Gen}}$ over $\mathfrak{B}_1$, universal with respect to $den^*_2$ and such that $den_2(T) \approx_{\mathfrak{B}_1, \mathfrak{B}_2} g_2$.

Piecha and Schroeder-Heister [Piecha \& Schroeder-Heister 2018] have recently disproved Prawitz's conjecture for valid arguments, although with respect to a notion of validity that differs from the one originally defined by Prawitz [Prawitz 1971b, 1973]. In our framework, this version corresponds to the following: $\Gamma \models^{\texttt{PS}} \alpha$ iff, for every base $\mathfrak{B}$, $\Gamma \models_\mathfrak{B} \alpha$.
\begin{thm}[Piecha \& Schroeder-Heister]
There are $\Gamma$ and $\alpha$ such that $\Gamma \models^{\texttt{PS}} \alpha$ and $\Gamma \nvdash_{\texttt{IL}} \alpha$.
\end{thm}
\noindent Piecha and Schroeder-Heister also suggest that there may hence be some intermediate logic $\mathfrak{L}$ which is $\texttt{PS}$-complete. We may thus ask if $\mathfrak{L}$ is also full-complete, as in conjecture 32. However, there may also be no such $\mathfrak{L}$, i.e. the class of rules which are $\texttt{PS}$-valid according to Prawitz's semantics could simply be non-recursive.

\subsection{Recognisability}

A crucial issue in Prawitz's proof-theoretic semantics is whether objects used to account for evidence are ‘‘epistemically transparent". In the semantics of valid arguments, for example, this is the problem of whether it is recognisable that a valid argument from - possibly empty - $\Gamma$ to $\alpha$ actually is a valid argument of this kind [see e.g. Prawitz 1973]. In ToG this is instead the problem of whether it is recognisable that a given construction actually is a ground for a (possibly general, hypothetical or general-hypothetical) assertion, and that a given term in a given language of grounding actually denotes such a ground. This problem is explicitly stated by Prawitz in many ground-theoretic papers [see mainly Prawitz 2015, 2019], and its centrality is easily seen as soon as one bears in mind what the fundamental task of ToG is: grounds cannot explain how valid inferences can be compelling if they are not recognisable, for no one would be compelled by an inference that simply yields an object of which we cannot also recognise that it justifies the conclusion - possibly under assumptions. We conclude this section with some observations about this topic, and about its relationship with the formalisation we have proposed in this paper.

It is not difficult to realise that, in ToG, the recognisability problem reduces to the problem of whether it is recognisable that an equation for an operation $f$ of type $\tau_1, ..., \tau_n \rhd \tau$ actually defines $f$ in such a way that, whenever applied to grounds for $\tau_i$ ($i \leq n$), $f$ yields a ground for $\tau$ - possibly after replacement of individual variables with individuals. Now, in Sections 4.2 and 4.3 we have said that operations on grounds must be \emph{total effective} functions. How should this be understood?

An operation on grounds $f$ of operational type $\tau_1, ..., \tau_n \rhd \tau$ can be said to be total if it is convergent on all values of its domain, i.e. if, for \emph{every} $g_i$ ground for $\tau_i$ ($i \leq n$), $f(g_1, ..., g_n)$ is a ground for $\tau$ - possibly after replacement of individual variables with individuals. The totality requirement is explicitly put forward by Prawitz [see mainly Prawitz 2015, 2019], and a main reason behind it may stem from the fact that operations on grounds are used, in ToG, to account for valid inferences. An operation on grounds associated to a valid inference \emph{must} be total, for a valid inference justifies us in asserting the conclusion \emph{whenever} we are justified in asserting the premises. If the operation is partial, there are cases in which we cannot infer the conclusion even if we have evidence for the premises; such an inference could not be considered as deductively correct, so that partial operations on grounds are of no use in ToG - or at least, we could well introduce them, but cannot also use them to analyse inferential validity.

However, the total character of operations on grounds is of course relevant when discussing the recognisability problem. If e.g. we require from the outset operations on grounds to be, in particular, total recursive functions, and restrict ourselves to composition methods which preserve totality, we obtain a decidable class of functions [see e.g. Dean 2021]. On the other hand, it is clear that the question whether a given equation defines a recursive function in such a way that this function is total cannot be expected to be decidable in general. What this ultimately shows is that we should now turn to our second question, i.e. what we mean by ‘‘effective".

An operation on grounds $f$ of operational type $\tau_1, ..., \tau_n \rhd \tau$ is expected to be effective in the sense of being effectively computable on all values of its domain, i.e. the defining equation must provide us with a set of instructions such that, for every $g_i$ ground for $\tau_i$ ($i \leq n$), we can effectively carry out the computation of $f(g_1, ..., g_n)$, thereby obtaining effectively a ground for $\tau$ - possibly after replacement of individual variables with individuals. This may suggest that - in line with the Church-Turing thesis - operations on grounds may be required to be recursive functions, or equivalent objects such as Turing machines, Church's $\lambda$-terms and so on. Is this requirement acceptable?

We remark that this suggestion is almost never contemplated by Prawitz. This may depend on a circularity highlighted by Peter with reference to the BHK-clauses [Peter 1959], and mentioned by Prawitz in the context of a discussion of a notion of proof which comes very close to that of ground [Prawitz 1977; for the similarity between Prawitz's BHK-like notion of proof and Prawitz's grounds see Tranchini 2014a, d'Aragona 2021b]. Peter's argument is that recursive functions cannot be used in a constructive determination of the meaning of the logical constants without falling into an inconsistency or into a circle. The definition of the notion of recursive function, in fact, must at some point employ an existential quantifier - e.g., existence of a scheme of equations that produces expected values, or existence of the zeros of the $\mu$-operator. This quantifier cannot be read in non-constructive terms, since the characterization of the meaning of the logical constants would in that case be inadequate; but if the existential must be read in necessarily constructive terms, the characterization would clearly be circular. From this Prawitz concludes that the notion of constructive function employed in the BHK-clauses for $\rightarrow$ and $\forall$ must be assumed to be primitive and not further analyzable. And this should also hold for the ground-clauses, as they too are taken to determine constructively the meaning of the logical constants.

Be that as it may, Prawitz seems to be sceptical about the possibility of obtaining a recognizability of the desired kind. With reference to terms of formal languages of grounding that denote (operations on) grounds, he highlights a distinction between \emph{knowing that} and \emph{knowing how}, which corresponds more or less to the distinction between canonical and non-canonical cases. In the canonical cases, an inferential agent

\begin{quote}
    know[s] that $T$ denotes a ground for asserting an atomic sentence $\alpha$ when this is how the meaning of $\alpha$ is given. Such knowledge is preserved by introduction inferences, given again that the meanings of the involved sentences are known: the term $T$ obtained by an introduction is in normal form, that is, it has the form $\phi(U)$ or $\phi(U, V)$, where $\phi$ is a primitive operation and the term $U$ or the terms $U$ and $V$ denote grounds for the premises - knowing that these terms do so, the agent also knows that $T$ denotes a ground for the conclusion. (Prawitz 2015, 97)
\end{quote}
So, the agent \emph{knows that} $T$ denotes a ground for $\alpha$ - provided she knows that the immediate sub-terms of $T$ denote grounds for the corresponding assertions, which may be granted by induction. By contrast, in the non-canonical cases

\begin{quote}
    the subject needs to reason from how $\phi$ is defined in order to see that $T$ denotes a ground for the conclusion. If $T$ is a closed term, she can in fact carry out the operations that $T$ is built up of and bring $T$ to normal form in this way, but she may not know this fact. Furthermore, when $T$ is an open term, it denotes a ground for an open assertion or an assertion under assumption, and it is first after appropriate substitutions for the free variables that one of the previous two cases arises. (Prawitz 2015, 97)
\end{quote}
Thus, the agent \emph{knows how} to obtain a ground for $\alpha$; indeed, $T$ amounts to a method that yields a ground for $\alpha$, possibly after application to individuals and to grounds for assumptions. But this knowledge-how may not also be a knowledge-that. In particular, we may assume that, when $T$ is a term in recursive language of grounding $\Lambda$, it is decidable whether $T$ denotes a ground for $\alpha$ - just as with recursive formal systems, where we can decide whether a given structure is a derivation in the system. But

\begin{quote}
    we know because of G\"{o}del’s incompleteness result that already for first order arithmetical assertions there is no closed language of grounds in which all grounds for them can be defined; for any such intuitively acceptable closed language of grounds, we can find an assertion and a ground for it that we find intuitively acceptable but that cannot be expressed within that language. (Prawitz 2015, 98)
\end{quote}
So languages of grounding must be understood as ‘‘open", i.e., as we said, as indefinitely expandable through addition of new operational symbols for new operations. And we seem to have no general argument granting that decidability \emph{within} $\Lambda$ is preserved throughout arbitrary expansions of $\Lambda$.

It is in any case clear that a thorough and more precise discussion of the recognisability problem requires that operations on grounds - and hence also languages of grounding and their expansions - be defined in a more rigorous way than what we have done in this paper. Some more refined characterisations may yield a class of operations on grounds defined by equations for which we can decide whether the operation they define is total. But these characterisations may turn out to be too restrictive given the aims of ToG. Other more refined characterisations may be sufficiently wide, but also imply undecidability.

As a concluding remark, we finally pinpoint that there may be independent reasons for renouncing the totality requirement on operations on grounds, for example in the context of a proof-theoretic treatment of paradoxes. A similar proposal has been for example put forward by some authors [see for example Tranchini 2019], although with reference to Prawitz's semantics of valid arguments.

\section{Conclusion}

The languages of grounding we have defined in this paper only involve terms. However, it is reasonable to ask of a formal language that it also contains formulas. A natural idea would be that of introducing two binary predicates - plus the logical constants: $Gr(T, \alpha)$ - meaning that $T$ denotes a ground for $\alpha$, or an operation with co-domain $\alpha$ if $T$ is open - and $T \equiv U$ - meaning that $T$ and $U$ denote equivalent (operations on) grounds. The approach in terms of denotation functions proposed here could be then further developed by giving a semantics for establishing whether formulas like $Gr(T, \alpha)$ or $T \equiv U$ hold, e.g. by giving truth- or proof-conditions for the formulas of the enriched language. In turn, this would require the introduction of a notion of model for a language of grounding. This line of thought may be explored in future works.

\medskip

\paragraph{Acknowledgments} I am grateful to Cesare Cozzo, Gabriella Crocco, Ansten Klev, Dag Prawitz and Peter Schroeder-Heister for helpful suggestions. I am also grateful to the anonymous reviewers, whose comments helped me to improve an earlier draft of this paper.

\begin{footnotesize}
\paragraph{References}
\begin{itemize}
\item C. Cozzo (1994), \emph{Meaning and argument. A theory of meaning centred on immediate argumental role}, Almqvist \& Wiksell, Uppsala.
\item A. P. d'Aragona (2018), \emph{A partial calculus for Dag Prawitz's theory of grounds and a decidability issue}, in in A. Christian, D. Hommen, N. Retzlaff, G. Schurz (eds), \emph{Philosophy of Science. European Studies in Philosophy of Science}, vol 9. Springer, Berlin Heidelberg New York.
\item --- (2019a), \emph{Dag Prawitz on proofs, operations and grounding}, in G. Crocco \& A. P. d'Aragona, \emph{Inferences and proofs}, special issue of \emph{Topoi}.
\item --- (2019b), \emph{Dag Prawitz's \emph{theory of grounds}}, PhD dissertation, Aix-Marseille University, ‘‘La Sapienza" University of Rome.
\item --- (2021a), \emph{Calculi of epistemic grounding based on Prawitz's theory of grounds}, submitted.
\item --- (2021b), \emph{Proofs, grounds and empty functions: epistemic compulsion in Prawitz's semantics}, in \emph{Journal of philosophical logic}, forthcoming.
\item W. Dean (2021), \emph{Recursive functions}, in E. N. Zalta (ed) \emph{The Stanford Encyclopedia of Philosophy} (Spring 2021 Edition).
\item G. F. Díez (2000), \emph{Five observation concerning the intended meaning of the intuitionistic logical constants}, in \emph{Journal of philosophical logic}, DOI : 10.1023/A:1004881914911.
\item K. Do\v{s}en (2015), \emph{Inferential semantics}, in H. Wansing (ed), \emph{Dag Prawitz on proofs and meaning}, Springer, Berlin Heidelberg New York.
\item M. Dummett (1991), \emph{The logical basis of metaphysics}, Harvard University Press, Cambridge.
\item --- (1993a), \emph{What is a theory of meaning \emph{(}I\emph{)}}, in M. Dummett, \emph{The seas of language}, Oxford University Press, Oxford.
\item --- (1993b), \emph{What is a theory of meaning \emph{(}II\emph{)}}, in M. Dummett, \emph{The seas of language}, Oxford University Press, Oxford.
\item N. Francez (2015), \emph{Proof-theoretic semantics}, College Publications, London.
\item G. Gentzen (1934 - 1935), \emph{Untersuchungen \"{u}ber das logische Schlie\ss{}en}, in \emph{Matematische Zeitschrift}, XXXIX, DOI: BF01201353.
\item A. Heyting(1956), \emph{Intuitionism. An introduction}, North-Holland Publishing Company, Amsterdam.
\item W. Howard (1980), \emph{The formula-as-types notion of construction}, in J. R. Hindley \& J. P. Seldin (eds), \emph{To H. B. Curry: essays on combinatoriy logic, lambda calculus and formalism}, Academic Press, London.
\item P. Martin-L\"{o}f (1984), \emph{Intuitionistic type theory}, Bibliopolis, Napoli.
\item R. Peter (1959), \emph{Rekursive Functionen}, Budapest, Akademiai Kiado.
\item T. Piecha (2016), \emph{Completeness in proof-theoretic semantics}, in T. Piecha \& P. Schroeder-Heister (eds), \emph{Advances in proof-theoretic semantics}, Springer, Berlin Heidelberg New York.
\item T. Piecha, W. de Campos Sanz \& P. Schroeder-Heister (2015), \emph{Failure of completeness in proof-theoretic semantics}, in \emph{Journal of philosophical logic}.
\item T. Piecha \& P. Schroeder-Heister (2018), \emph{Incompleteness of intuitionistic propositional logic with respect to proof-theoretic semantics}, in \emph{Studia Logica}.
\item F. Poggiolesi (2016), \emph{A critical overview of the most recent logics of grounding}, in F. Boccuni \& A. Sereni (eds), \emph{Objectivity, realism and proof}, Boston Studies in the Philosophy and History of Science, Springer, 2016.
\item --- (2020), \emph{Logics of grounding}, in M. Raven (ed), \emph{Routledge handbook for metaphysical grounding}, pp. 213 - 227, Routledge, 2020.
\item D. Prawitz (1971a), \emph{Constructive semantics}, in \emph{Proceedings of the First
Scandinavian Logic Symposium}.
\item --- (1971b), \emph{Ideas and results in proof theory}, in \emph{Proceedings of the Second Scandinavian Logic Symposium}.
\item --- (1973), \emph{Towards a foundation of a general proof-theory}, in P. Suppes, (ed) \emph{Logic methodology and philosophy of science IV}, North-Holland Publishing Company, Amsterdam, DOI: S0049-237X(09)70361-1.
\item --- (1977), \emph{Meaning and proofs: on the conflict between classical and intuitionistic logic}, in \emph{Theoria}.
\item --- (1979), \emph{Proofs and the meaning and completeness of the logical constants}, in J. Hintikka et al (eds), \emph{Essays on mathematical and philosophical logic}, Reidel, Dordrecht.
\item --- (1977), \emph{Meaning and proofs: on the conflict between classical and intuitionistic logic}, in \emph{Theoria}, DOI: j.1755-2567.1977.tb00776.x.
\item --- (2006), \emph{Natural deduction. A proof-theoretical study}, Dover, New York.
\item --- (2009), \emph{Inference and knowledge}, in M. Pelis (ed), \emph{The Logica Yearbook 2008}, College Publications, London.
\item --- (2012a), \emph{The epistemic significance of valid inference}, in \emph{Synthese}, DOI: s11229-011-9907-7.
\item --- (2012b), \emph{Truth and proof in intuitionism}, in P. Dybier, S. Lindström, E. Palmgren \& G. Sundholm (eds), \emph{Epistemology versus ontology}, Springer, Berlin Heidelberg New York.
\item --- (2013), \emph{Validity of inferences}, in M. Frauchiger (ed), \emph{Reference, rationality, and phenomenology: themes from F\o{}llesdal}, Dordrecht, Ontos Verlag.
\item --- (2014), \emph{An approach to general proof theory and a conjecture of a kind of completeness of intuitionistic logic revisited}, in L.C. Pereira, E. Heusler and V. de Paiva (eds), \emph{Advances in
natural deduction}, Trends in Logic (Studia Logica Library), vol. 39, Springer, Dordrecht, 2014.
\item --- (2015), \emph{Explaining deductive inference}, in H. Wansing (ed), \emph{Dag Prawitz on proofs and meaning}, Springer, Berlin Heidelberg New York, DOI: 978-3-319-11041-7\_3.
\item --- (2019), \emph{The seeming interdependence between the concepts of valid inference and proof}, in G. Crocco \& A. Piccolomini d'Aragona (eds), \emph{Inferences and proofs}, special issue of \emph{Topoi}, DOI: s11245-017-9506-4.
\item --- (2020a), \emph{Validity of inferences}, forthcoming.
\item --- (2020b), \emph{Validity of inferences reconsidered}, forthcoming.
\item P. Schroeder-Heister (1984a), \emph{A natural extension for natural deduction}, in \emph{Journal of symbolic logic}.
\item --- (1984b), \emph{Generalized rules for quantifiers and the completeness of the intuitionistic operators $\wedge$, $\vee$, $\rightarrow$, $\forall$, $\exists$}, in M. M. Richter, E. B\"{o}rger, W. Oberschelp, B. Schinzel \& W. Thomas (eds.), \emph{Computation and proof theory. Proceedings of the Logic Colloquium held in Aachen, July 18-23, 1983, Part II}, Springer, Berlin Heidelberg New York Tokyo.
\item --- (1991), \emph{Uniform proof-theoretic semantics for logical constants. Abstract}, in \emph{Journal of Symbolic Logic}.
\item --- (2006), \emph{Validity concepts in proof-theoretic semantics}, in \emph{Synthese}.
\item --- (2008), \emph{Proof-theoretic versus model-theoretic consequence}, in M. Peli\v{s} (ed), \emph{The Logica Yearbook 2008}, Filosofia, Prague.
\item --- (2012), \emph{The categorical and the hypothetical: a critique of some fundamental assumptions of standard semantics}, in S. Lindstr\"{o}m et al. (eds), \emph{The philosophy of logical consequence and inference. Proceedings of the workshop "The Philosophy of Logical Consequence, Uppsala, 2008"}, Synthese.
\item --- (2018), \emph{Proof-theoretic semantics}, in Edward N. Zalta (ed.), \emph{The Stanford Encyclopedia of Philosophy \emph{(}Spring 2018 Edition\emph{)}}.
\item G. Sundholm (1998), \emph{Proofs as acts and proofs as objects}, in \emph{Theoria}.
\item L. Tranchini (2014a), \emph{Dag Prawitz}, APhEx 9.
\item --- (2014b), \emph{Proof-theoretic semantics, proofs and the distinction between sense and denotation}, in \emph{Journal of logic and computation}.
\item --- (2019), \emph{Proof, meaning and paradox. Some remarks}, in G. Crocco \& A. Piccolomini d’Aragona (eds), \emph{Inferences and proofs, special issue of Topoi}.
\item A. S. Troelstra \& D. Van Dalen (1988), \emph{Constructivism in mathematics}, vol. I, North-Holland Publishing Company, Amsterdam.
\item G. Usberti (2015), \emph{A notion of $C$-justification for empirical statements}, in H. Wansing (ed), \emph{Dag Prawitz on proofs and meaning}, Springer, Berlin Heidelberg New York, DOI: 978-3-319-11041-7\_18.
\item --- (2019), \emph{Inference and epistemic transparency}, in G. Crocco \& A. Piccolomini d'Aragona (eds), \emph{Inferences and proofs}, special issue of \emph{Topoi}, DOI: s11245-017-9497-1.
\end{itemize}
\end{footnotesize}
\end{document}